\documentclass[11pt]{article}
\usepackage{amsfonts}

\usepackage[dvips]{graphicx}
\usepackage{amsmath}

\newtheorem{theorem}{Theorem}[subsection]

\newtheorem{conjecture}[theorem]{Conjecture}
\newtheorem{corollary}[theorem]{Corollary}

\newtheorem{lemma}[theorem]{Lemma}

\newtheorem{proposition}[theorem]{Proposition}

\newenvironment{proof}[1][Proof]{\textbf{#1.} }{\ \rule{0.5em}{0.5em}}

\input{tcilatex}

\begin{document}

\title{Branching rules, Kostka-Foulkes polynomials and $q$-multiplicities in tensor
product for the root systems $B_{n},C_{n}$ and $D_{n}$}
\author{C\'{e}dric Lecouvey \\
%EndAName
lecouvey@math.unicaen.fr}
\date{}
\maketitle

\begin{abstract}
The Kostka-Foulkes polynomials $K_{\lambda,\mu}^{\phi}(q)$ related to a root
system $\phi$ can be defined as alternated sums running over the Weyl group
associated to $\phi.$ By restricting these sums over the elements of the
symmetric group when $\phi$ is of type $B_{n},C_{n}$ or $D_{n}$, we obtain
again a class $\widetilde{K}_{\lambda,\mu}^{\phi}(q)$ of Kostka-Foulkes
polynomials. When $\phi$ is of type $C_{n}$ or $D_{n}$ there exists a
duality beetween these polynomials and some natural $q$-multiplicities $%
u_{\lambda ,\mu}(q)$ and $U_{\lambda,\mu}(q)$ in tensor product \cite{lec}.
In this paper we first establish identities for the $\widetilde{K}%
_{\lambda,\mu}^{\phi}(q)$ which implies in particular that they can be
decomposed as sums of Kostka-Foulkes polynomials $K_{\lambda,%
\mu}^{A_{n-1}}(q)$ with nonnegative integer coefficients.\ Moreover these
coefficients are branching rule coefficients$.$ This allows us to clarify
the connection beetween the $q$-multiplicities $u_{\lambda,\mu}(q),U_{%
\lambda,\mu}(q)$ and the polynomials $K_{\lambda,\mu}^{\diamondsuit}(q)$
defined in \cite{SZ}. Finally we show that $u_{\lambda,\mu}(q)$ and $%
U_{\lambda,\mu}(q)$ coincide up to a power of $q$ with \ the one dimension
sum introduced in \cite{Ok} when all the parts of $\mu$ are equal to $1$
which partially proves some conjectures of \cite{lec} and \cite{SZ}.
\end{abstract}

\section{Introduction}

Consider $\lambda$ and $\mu$ two partitions of the set $\mathcal{P}_{n}$ of
partitions with $n$ parts.\ The Schur-Weyl duality establishes that the
dimension $K_{\lambda,\mu}^{A_{n-1}}$ of the weight space $\mu$ in the
finite dimensional irreducible $sl_{n}$-module $V^{A_{n-1}}(\lambda)$ of
highest weight $\lambda$ is equal to the multiplicity of $%
V^{A_{n-1}}(\lambda)$ in the tensor product 
\begin{equation*}
V_{\mu}^{A_{n-1}}=V^{A_{n-1}}(\mu_{1}\Lambda_{1})\otimes\cdot\cdot\cdot%
\otimes V^{A_{n-1}}(\mu_{n}\Lambda_{1}).
\end{equation*}
It follows from the Weyl character formula that $K_{\lambda,\mu}^{A_{n-1}}=%
\sum_{\sigma\in\mathcal{S}_{n}}(-1)^{l(\sigma)}\mathcal{P}%
^{A_{n-1}}(\sigma(\lambda+\rho)-(\mu+\rho))$ where $\mathcal{P}^{A_{n-1}}$
is the Kostant partition function which counts, in the root system of type $%
A_{n-1},$ the number of decomposition of $\beta\in\mathbb{Z}^{n}$ as a sum
of positive roots. The Kostka-Foulkes polynomials can be defined by setting $%
K_{\lambda ,\mu}^{A_{n-1}}(q)=\sum_{\sigma\in\mathcal{S}_{n}}(-1)^{l(\sigma)}%
\mathcal{P}_{q}^{A_{n}}(\sigma(\lambda+\rho)-(\mu+\rho))$ where $\mathcal{P}%
_{q}^{A_{n-1}}$ is the $q$-Kostant partition function characterized by 
\begin{equation*}
\prod_{\alpha\text{ positive root}}\dfrac{1}{(1-qe^{\alpha})}=\sum_{\beta \in%
\mathbb{Z}^{n}}\mathcal{P}_{q}^{A_{n}}(\beta)e^{\beta}
\end{equation*}
with $\rho=(n-1,...,0)$ the half sum of the positive roots.\ One can prove
that they are the coefficients of the expansion of the Schur function $%
s_{\mu }(x)$ on the basis of Hall polynomials $\{P_{\lambda}(x,q),\lambda \in%
\mathcal{P}_{n}\}$ (see \cite{mac})$.\;$Then it follows from the theory of
affine Hecke algebras that the Kostka-Foulkes polynomials are
Kazhdan-Lusztig polynomials \cite{Lu}.\ In particular they have nonnegative
integer coefficients. As proved by Lascoux and Sch\"{u}tzenberger this
positivity result can also be obtained by using the charge statistic $%
\mathrm{ch}$ on semistandard tableaux. More precisely we have $%
K_{\lambda,\mu}^{A_{n-1}}(q)=\sum_{T\in ST(\lambda)_{\mu}}q^{\mathrm{ch}(T)}$
where $ST(\lambda)_{\mu }$ is the set of semistandard tableaux of shape $%
\lambda$ and weight $\mu.$ In \cite{NY}, Nakayashiki and Yamada have shown
that the charge can be computed from the combinatorial $R$-matrix
corresponding to Kashiwara's crystals associated to some $U_{q}(\widehat{%
sl_{n}})$-modules.

Now consider $\phi$ a root system of type $B_{n},C_{n}$ or $D_{n}$.\ Write $%
\frak{g}_{\phi}$ for the corresponding simple Lie algebra.\ The
Kostka-Foulkes polynomials $K_{\lambda,\mu}^{\phi}(q)$ associated to $\phi$
are defined by setting 
\begin{equation*}
K_{\lambda,\mu}^{\phi}(q)=\sum_{w\in W_{\phi}}(-1)^{l(w)}\mathcal{P}%
_{q}^{\phi}(w(\lambda+\rho_{\phi})-(\mu+\rho_{\phi}))
\end{equation*}
where $W_{\phi},\rho_{\phi}$ and $\mathcal{P}_{q}^{\phi}$ are respectively
the Weyl group, the half sum of the positive roots and the $q$-partition
function corresponding to $\phi$. The polynomial $K_{\lambda,\mu}^{\phi}(q)$
can be considered as a $q$-analogue of the dimension of the weight space $%
\mu $ in $V^{\phi}(\lambda).$ As Kazhdan-Lusztig polynomials, they have also
nonnegative coefficients. In \cite{Lec3}, we have obtained for the root
systems $B_{n},C_{n}$ and $D_{n}$ a statistic on Kashiwara-Nakashima's
tableaux from which it is possible to deduce this positivity for particular
pairs of partitions $(\lambda,\mu)$. Nevertheless as far as the author is
aware, no combinatorial proof of this positivity result is known in general.

\noindent The Kostka-Foulkes polynomials $K_{\lambda ,\mu }^{\phi }(q)$
cannot be directly interpreted as $q$-multiplicities in tensor products. So
there does not exist an equivalent result to the Schur Weyl duality for the
root system $\phi .\;$Denote by $V^{\phi }(\lambda )$ the finite dimensional
irreducible $\frak{g}_{\phi }$-module of highest weight $\lambda .\;$In \cite
{lec}, we have introduced from determinantal expressions of the Schur
functions associated to $\phi $, two polynomials $u_{\lambda ,\mu }(q)$ and $%
U_{\lambda ,\mu }(q)$ which can be respectively regarded as quantizations of
the multiplicities of $V^{\phi }(\lambda )$ in the tensor products 
\begin{equation*}
V_{\mu }^{\phi }=V^{\phi }(\mu _{1}\Lambda _{1})\otimes \cdot \cdot \cdot
\otimes V^{\phi }(\mu _{n}\Lambda _{1})\text{ and }W_{\mu }^{\phi }=W^{\phi
}(\mu _{1}\Lambda _{1})\otimes \cdot \cdot \cdot \otimes W^{\phi }(\mu
_{n}\Lambda _{1})
\end{equation*}
where for any $i=1,...,n,$ $W^{\phi }(\mu _{i}\Lambda _{1})=V^{\phi }(\mu
_{i}\Lambda _{1})\oplus V^{\phi }((\mu _{i}-2)\Lambda _{1})\oplus \cdot
\cdot \cdot \oplus V^{\phi }((\mu _{i}\mathrm{mod}2)\Lambda _{1})$.\ When $n$
is sufficiently large they do not depend on the root system $\phi $
considered and we have established a duality result between the $q$%
-multiplicities $u_{\lambda ,\mu }(q),$ $U_{\lambda ,\mu }(q)$ and the
polynomials 
\begin{equation*}
\widetilde{K}_{\lambda ,\mu }^{\phi }(q)=\sum_{\sigma \in \mathcal{S}%
_{n}}(-1)^{l(\sigma )}\mathcal{P}_{q}^{\phi }(w(\lambda +\rho _{n})-(\mu
+\rho _{n}))
\end{equation*}
where $\rho _{n}=(n,...,1).\;$These polynomials $\widetilde{K}_{\lambda ,\mu
}^{\phi }(q)$ are also Kostka-Foulkes polynomials. So this result can be
interpreted as a duality between $q$-analogues of weight multiplicities and $%
q$-analogues of tensor product multiplicities for the root systems $%
B_{n},C_{n}$ and $D_{n}$.

\noindent At the same time Shimozono and Zabrocki \cite{SZ} have
independently defined by using creating operators some polynomials $%
K_{\lambda ,R}^{\diamondsuit }(q)$ where $R$ is a sequence of rectangular
partitions and $\diamondsuit $ a partition of the set $\{\emptyset
,(1),(11),(2)\}$.\ These polynomials can also be regarded as $q$%
-multiplicities in tensor products. In \cite{Ok}, Hatayama, Kuniba, Okado
and Takagi have introduced for type $C_{n}$ a quantization $X_{\lambda ,\mu
}(q)$ of the multiplicity of $V^{C_{n}}(\lambda )$ in $W_{\mu }^{C_{n}}.\;$%
This quantization is based on the determination of the combinatorial $R$%
-matrix of some $U_{q}^{\prime }(\widehat{sp_{2n}})$-crystals in the spirit
of \cite{NY}. It can be regarded as a one dimension sum for the affine root
system $C_{n}^{(1)}.$ In \cite{lec} and \cite{SZ}, the authors conjecture
that the polynomials $X_{\lambda ,\mu }(q),$ $U_{\lambda ,\mu }(q)$ and $%
K_{\lambda ,\mu }^{(2)}(q)$ coincide up to simple renormalizations$.$ As
observed in \cite{SZ}, this conjecture can be related to the $X=M$
conjecture which gives fermionic formulas for the one dimension sum $X$.\
Note that the $X=M$ conjecture have been proved in various cases for all
nonexceptional affine types \cite{OSS}, \cite{ok3} and \cite{SCS}$.$

In this article we first obtained identities for the polynomials $\widetilde{%
K}_{\lambda ,\mu }^{\phi }(q)$ which imply that they can be decomposed as
sums of polynomials $K_{\lambda ,\mu }^{A_{n-1}}(q)$.\ Moreover the
coefficients of these decompositions can be simply expressed in terms of
branching rules coefficients. This gives in particular an elementary proof
of the positivity of the Kostka-Foulkes polynomials $\widetilde{K}_{\lambda
,\mu }^{\phi }(q).$ Next we obtain similar decompositions for the
polynomials $u_{\lambda ,\mu }(q)$ and $U_{\lambda ,\mu }(q).$ By comparing
these identities with those obtained for the polynomials $K_{\lambda
,R}^{\diamondsuit }(q)$ in \cite{SZ}, we derive the equalities $K_{\lambda
,\mu }^{(1,1)}(q)=u_{\lambda ,\mu }(q^{2})$ and $K_{\lambda ,\mu
}^{(2)}(q)=U_{\lambda ,\mu }(q^{2}).$ Finally we establish some conjectures
of \cite{lec} and \cite{SZ} when all the parts of $\mu $ are equal to $1$
(i.e. for the $q$-multiplicities defined in the tensor powers of the vector
representation), namely we have 
\begin{equation}
K_{\lambda ,(1^{n})}^{(1,1)}(q)=u_{\lambda ,(1^{n})}(q^{2})=q^{n-\left|
\lambda \right| }X_{\lambda ,(1^{n})}(q^{2})\text{ and }K_{\lambda
,(1^{n})}^{(2)}(q)=U_{\lambda ,(1^{n})}(q^{2})=q^{2(n-\left| \lambda \right|
)}X_{\lambda ,(1^{n})}(q^{2}).  \label{EQF}
\end{equation}

\bigskip

In Section $2$ we review some material on root systems, branching rules
coefficients, Kostka-Foulkes polynomials$\ $and $q$-multiplicities $%
u_{\lambda ,\mu }(q),U_{\lambda ,\mu }(q)$ we need in the sequel. In section 
$3$ we obtain identities for the polynomials $\widetilde{K}_{\lambda ,\mu
}^{\phi }(q),$ $u_{\lambda ,\mu }(q)$ and $U_{\lambda ,\mu }(q)$ from which
we clarify the relations between $u_{\lambda ,\mu }(q),U_{\lambda ,\mu }(q)$
and $K_{\lambda ,\mu }^{(1,1)}(q),K_{\lambda ,\mu }^{(2)}(q).$ Section $4$
is devoted to the proof of (\ref{EQF}). Note that the $X=M$ conjecture is in
particular true when all the parts of $\mu $ are equal to $1$ \cite{OSS}.\
Thus in this case the one dimension sums $X$ and their corresponding
fermionic formulas $M$ are, up to simple renormalizations, Kazhdan-Lusztig
polynomials.

\bigskip

\noindent\textbf{Notation: }In the sequel we frequently define similar
objects for the root systems $B_{n}$ $C_{n}$ and $D_{n}$. When they are
related to type $B_{n}$ (resp. $C_{n},D_{n}$), we implicitly attach to them
the label $B$ (resp. the labels $C,D$). To avoid cumbersome repetitions, we
sometimes omit the labels $B,C$ and $D$ when our definitions or statements
are identical for the three root systems.

\section{Background}

\subsection{Convention for the root systems of types $B_{n},C_{n}$ and $%
D_{n} $}

Consider an integer $n\geq1.$ The weight lattice for the root system $C_{n}$
(resp. $B_{n}$ and $D_{n})$ can be identified with $P_{C_{n}}=\mathbb{Z}^{n}$
(resp. $P_{B_{n}}=P_{D_{n}}=\left( \dfrac{\mathbb{Z}}{2}\right) ^{n})$
equipped with the orthonormal basis $\varepsilon_{i},$ $i=1,...,n$.\ We take
for the simple roots 
\begin{equation}
\left\{ 
\begin{tabular}{l}
$\alpha_{n}^{B_{n}}=\varepsilon_{n}\text{ and }\alpha_{i}^{B_{n}}=%
\varepsilon_{i}-\varepsilon_{i+1}\text{, }i=1,...,n-1\text{ for the root
system }B_{n}$ \\ 
$\alpha_{n}^{C_{n}}=2\varepsilon_{n}\text{ and }\alpha_{i}^{C_{n}}=%
\varepsilon_{i}-\varepsilon_{i+1}\text{, }i=1,...,n-1\text{ for the root
system }C_{n}$ \\ 
$\alpha_{n}^{D_{n}}=\varepsilon_{n}+\varepsilon_{n-1}\text{ and }\alpha
_{i}^{D_{n}}=\varepsilon_{i}-\varepsilon_{i+1}\text{, }i=1,...,n-1\text{ for
the root system }D_{n}$%
\end{tabular}
\right. .  \label{simple_roots}
\end{equation}
Then the set of positive roots are 
\begin{equation*}
\left\{ 
\begin{tabular}{l}
$R_{B_{n}}^{+}=\{\varepsilon_{i}-\varepsilon_{j},\varepsilon_{i}+%
\varepsilon_{j}\text{ with }1\leq i<j\leq n\}\cup\{\varepsilon_{i}\text{
with }1\leq i\leq n\}\text{ for the root system }B_{n}$ \\ 
$R_{C_{n}}^{+}=\{\varepsilon_{i}-\varepsilon_{j},\varepsilon_{i}+%
\varepsilon_{j}\text{ with }1\leq i<j\leq n\}\cup\{2\varepsilon_{i}\text{
with }1\leq i\leq n\}\text{ for the root system }C_{n}$ \\ 
$R_{D_{n}}^{+}=\{\varepsilon_{i}-\varepsilon_{j},\varepsilon_{i}+%
\varepsilon_{j}\text{ with }1\leq i<j\leq n\}\text{ for the root system }%
D_{n}$%
\end{tabular}
\right. .
\end{equation*}
Denote respectively by $P_{B_{n}}^{+},P_{C_{n}}^{+}$ and $P_{D_{n}}^{+}$the
sets of dominant weights of $so_{2n+1},sp_{2n}$ and $so_{2n}.$

\noindent Let $\lambda=(\lambda_{1},...,\lambda_{n})$ be a partition with $n$
parts. We will classically identify $\lambda$ with the dominant weight $%
\sum_{i=1}^{n}\lambda_{i}\varepsilon_{i}.$ Note that there exists dominant
weights associated to the orthogonal root systems whose coordinates on the
basis $\varepsilon_{i},$ $i=1,...,n$ are not positive integers (hence which
cannot be regarded as partitions). For each root system of type $B_{n},C_{n}$
or $D_{n},$ the set of weights having nonnegative integer coordinates on the
basis $\varepsilon_{1},...,\varepsilon_{n}$ can be identify with the set $%
\mathcal{P}_{n}$ of partitions of length $n.$\ For any partition $\lambda,$
the weights of the finite dimensional $so_{2n+1},sp_{2n}$ or $so_{2n}$%
-module of highest weight $\lambda$ are all in $\mathbb{Z}^{n}.\;$For any $%
\alpha \in\mathbb{Z}^{n}$ we write $\left| \alpha\right| =\alpha_{1}+\cdot
\cdot\cdot+\alpha_{n}$ and $\left\| \alpha\right\|
=\sum_{i=1}^{n-1}(n-i)\alpha_{i}.$

\noindent The conjugate partition of the partition $\lambda$ is denoted $%
\lambda^{\prime}$ as usual.\ Consider $\lambda,\mu$ two partitions of length 
$n$ and set $m=\max(\lambda_{1},\mu_{1})$.\ Then by adding to $\lambda
^{\prime}$ and $\mu^{\prime}$ the required numbers of parts $0$ we will
consider them as partitions of length $m.$

\bigskip

\noindent The Weyl group $W_{B_{n}}=W_{C_{n}}$ of $so_{2n+1}$ and $sp_{2n}$
is identified to the subgroup of the permutation group of the set $%
\{\overline {n},...,\overline{2},\overline{1},1,2,...,n\}$\ generated by $%
s_{i}=(i,i+1)(\overline{i},\overline{i+1}),$ $i=1,...,n-1$ and $s_{n}=(n,%
\overline{n})$ where for $a\neq b$ $(a,b)$ is the simple transposition which
switches $a$ and $b.$ We denote by $l_{B}$ the length function corresponding
to the set of generators $s_{i},$ $i=1,...n.$

\noindent The Weyl group $W_{D_{n}}$ of $so_{2n}$ is identified to the
subgroup of $W_{B_{n}}$\ generated by the transpositions $s_{i}=(i,i+1)(%
\overline{i},\overline{i+1}),$ $i=1,...,n-1$ and $s_{n}^{\prime }=(n,%
\overline{n-1})(n-1,\overline{n})$. We denote by $l_{D}$ the length function
corresponding to the set of generators $s_{n}^{\prime}$ and $s_{i},$ $%
i=1,...n-1.$

\noindent Note that $W_{D_{n}}\subset W_{B_{n}}$ and any $w\in W_{B_{n}}$
verifies $w(\overline{i})=\overline{w(i)}$ for $i\in\{1,...,n\}.$ The action
of $w$ on $\beta=(\beta_{1},...,\beta_{n})\in\mathbb{Z}^{n}$ is given by 
\begin{equation*}
w\cdot(\beta_{1},...,\beta_{n})=(\beta_{1}^{w},...,\beta_{n}^{w})
\end{equation*}
where $\beta_{i}^{w}=\beta_{w(i)}$ if $\sigma(i)\in\{1,...,n\}$ and $\beta
_{i}^{w}=-\beta_{w(\overline{i})}$ otherwise.

\noindent The half sums $\rho_{B_{n}},\rho_{C_{n}}$ and $\rho_{D_{n}}$ of
the positive roots associated to each root system $B_{n},C_{n}$ and $D_{n}$
verify: 
\begin{equation*}
\rho_{B_{n}}=(n-\dfrac{1}{2},n-\dfrac{3}{2},...,\dfrac{1}{2}%
),\rho_{C_{n}}=(n,n-1,...,1)\text{ and }\rho_{B_{n}}=(n-1,n-2,...,0).
\end{equation*}
In the sequel we identify the symmetric group $\mathcal{S}_{n}$ (which is
the Weyl group of the root system $A_{n-1})$ with the subgroup of $W_{B_{n}}$
or $W_{D_{n}}$ generated by the $s_{i}$'s, $i=1,...,n-1.$

\subsection{Branching rules coefficients}

For any partition $\lambda,$ we denote by $V_{n}^{B}(\lambda),V_{n}^{C}(%
\lambda),$ and $V_{n}^{D}(\lambda)$ the finite dimensional irreducible
modules of highest weight $\lambda$ respectively for $sp_{2n},so_{2n+1}$ and 
$so_{2n}.$ Then $V_{n}^{B}(\lambda),V_{n}^{C}(\lambda),$ and $%
V_{n}^{D}(\lambda)$ can also be regarded as irreducible representations
respectively of the groups $Sp_{2n},So_{2n+1}$ and $So_{2n}.$ By restriction
to $GL_{n},$ they decompose in a direct sum of irreducible rational
representations. Recall that the irreducible rational representations of $%
GL_{n}$ are indexed by the $n$-tuples 
\begin{equation}
(\gamma^{+},\gamma^{-})=(\gamma_{1}^{+},\gamma_{2}^{+},...,%
\gamma_{p}^{+},0,...,0,-\gamma_{q}^{-},...,-\gamma_{1}^{-})  \label{jamma+-}
\end{equation}
where $\gamma^{+}$ and $\gamma^{-}$ are partitions of length $p$ and $q$
such that $p+q\leq n.$ Write $V_{n}^{A}(\gamma^{+},\gamma^{-})$ for the
irreducible rational representations of $GL_{n}$ of highest weight $%
(\gamma^{+},\gamma ^{-}).$ When $\gamma^{-}=\emptyset,$ we write simply $%
V_{n}^{A}(\gamma)$ instead of $V_{n}^{A}(\gamma^{+},\gamma^{-}).$

\bigskip

As customary, we use for a basis of the group algebra $\mathbb{Z}[\mathbb{Z}%
^{n}],$ the formal exponentials $(e^{\beta})_{\beta\in \mathbb{Z}^{n}}$
satisfying the relations $e^{\beta_{1}}e^{\beta_{2}}=e^{\beta_{1}+%
\beta_{2}}. $ We furthermore introduce $n$ independent indeterminates $%
x_{1},...,x_{n}$ in order to identify $\mathbb{Z}[\mathbb{Z}^{n}]$ with the
ring of polynomials $\mathbb{Z}[x_{1},...,x_{n},x_{1}^{-1},...,x_{n}^{-1}]$
by writing $e^{\beta}=x_{1}^{\beta_{1}}\cdot\cdot\cdot
x_{n}^{\beta_{n}}=x^{\beta}$ for any $\beta=(\beta _{1},...,\beta_{n})\in%
\mathbb{Z}^{n}.$

\noindent Set 
\begin{gather}
\prod_{1\leq r<s\leq n}\left( 1-\frac{1}{x_{r}x_{s}}^{-1}\right)
\prod_{1\leq i\leq n}\left( 1-\dfrac{1}{x_{i}}\right) ^{-1}=\sum_{\beta\in
L_{B}}b(\beta)x^{-\beta},  \label{def_b} \\
\prod_{1\leq r\leq s\leq n}\left( 1-\frac{1}{x_{r}x_{s}}\right)
^{-1}=\sum_{\beta\in L_{C}}c(\beta)x^{-\beta}\text{,}  \notag \\
\prod_{1\leq r<s\leq n}\left( 1-\frac{1}{x_{r}x_{s}}\right)
^{-1}=\sum_{\beta\in L_{D}}d(\beta)x^{-\beta}  \notag
\end{gather}
where 
\begin{gather*}
L_{B}=\{\beta\in\mathbb{Z}^{n},\beta=\sum_{1\leq r<s\leq
n}e_{r,s}(\varepsilon_{r}+\varepsilon_{s})+\sum_{1\leq i\leq
n}e_{i}\varepsilon _{i}\text{ with }e_{r,s}\geq0\text{ and }e_{i}\geq0\} \\
L_{C}=\{\beta\in\mathbb{Z}^{n},\beta=\sum_{1\leq r\leq s\leq
n}e_{r,s}(\varepsilon_{r}+\varepsilon_{s})\text{ with }e_{r,s}\geq0\}\text{
and } \\
L_{D}=\{\beta\in\mathbb{Z}^{n},\beta=\sum_{1\leq r<s\leq
n}e_{r,s}(\varepsilon_{r}+\varepsilon_{s})\text{ with }e_{r,s}\geq0\}.
\end{gather*}
Denote respectively by $[V_{n}^{A}(\gamma^{+},\gamma^{-}):V_{n}^{B}(%
\lambda)],$ $[V_{n}^{A}(\gamma^{+},\gamma^{-}):V_{n}^{C}(\lambda)]$ and $%
[V_{n}^{A}(\gamma^{+},\gamma^{-}):V_{n}^{D}(\lambda)],$ the multiplicities
of $V_{n}^{A}(\gamma^{+},\gamma^{-})$ in the restrictions of $%
V_{n}^{B}(\lambda),V_{n}^{C}(\lambda)$ and $V_{n}^{D}(\lambda)$ to $GL_{n}.$

\begin{proposition}
\label{prop_mul_sum}With the above notation, we have:

\begin{enumerate}
\item  $[V_{n}^{A}(\gamma ^{+},\gamma ^{-}):V_{n}^{B}(\lambda )]=\sum_{w\in
W_{B_{n}}}(-1)^{l(w)}b(w\circ \lambda -(\gamma ^{+},\gamma ^{-})),$

\item  $[V_{n}^{A}(\gamma ^{+},\gamma ^{-}):V_{n}^{C}(\lambda )]=\sum_{w\in
W_{C_{n}}}(-1)^{l(w)}c(w\circ \lambda -(\gamma ^{+},\gamma ^{-})),$

\item  $[V_{n}^{A}(\gamma ^{+},\gamma ^{-}):V_{n}^{D}(\lambda )]=\sum_{w\in
W_{D_{n}}}(-1)^{l(w)}d(w\circ \lambda -(\gamma ^{+},\gamma ^{-})).$
\end{enumerate}
\end{proposition}

\begin{proof}
The proposition can be considered as a corollary of Theorem 8.2.1 of \cite
{GW} with $G$ one of the Lie groups $So_{2n+1},Sp_{2n},So_{2n}$ and $%
H=GL_{n}.$
\end{proof}

\bigskip

\noindent For any partitions $\lambda$ and $\nu$ of length $n,$ write $%
[V_{n}^{D}(\nu):V_{n}^{B}(\lambda)]$ for the multiplicity of $V_{n}^{D}(\nu)$
in the restriction of $V_{n}^{B}(\lambda)$ to $So_{2n}.$

\begin{lemma}
\label{lem_dec_b_D}With the notation above we have 
\begin{equation*}
\lbrack V_{n}^{D}(\nu ):V_{n}^{B}(\lambda )]=\sum_{w\in W_{B_{n}}}(-1)^{l(w)}%
\mathbf{1}_{\mathbb{N}}(w\circ \lambda -\nu )
\end{equation*}
where for any $\beta \in \mathbb{Z}^{n},$ $\mathbf{1}_{\mathbb{N}}(\beta )=1$
if all the coordinates of $\beta $ are nonnegative integers and $\mathbf{1}_{%
\mathbb{N}}(\beta )=0$ otherwise.
\end{lemma}

\begin{proof}
The lemma also follows from Theorem 8.2.1 of \cite{GW}.
\end{proof}

\noindent For any partitions $\lambda$ and $\nu$ of length $n,$ write $%
[V_{n}^{B}(\lambda):V_{2n}^{A}(\nu)]$ for the multiplicity of $%
V_{n}^{B}(\lambda)$ in the restriction of $V_{2n}^{A}(\nu)$ from $GL_{2n}$
to $So_{2n}.$ Similarly write $[V_{n}^{C}(\lambda):V_{2n}^{A}(\nu)]$ for the
multiplicity of $V_{n}^{C}(\lambda)$ in the restriction of $V_{2n}^{A}(\nu)$
from $GL_{2n}$ to $Sp_{2n}$ and $[V_{n}^{D}(\lambda):V_{2n}^{A}(\nu)]$ for
the multiplicity of $V_{n}^{D}(\lambda)$ in the restriction of $%
V_{2n}^{A}(\nu)$ from $GL_{2n}$ to $So_{2n+1}.$ Denote respectively by $%
\mathcal{P}_{n}^{(2)}$ and $\mathcal{P}_{n}^{(1,1)}$ the sub-sets of $%
\mathcal{P}_{n}$ containing the partitions with even rows and the partitions
with even columns. The Littlewood-Richardson coefficients are denoted $%
c_{\gamma,\lambda}^{v}$ as usual.

\noindent Let us recall a classical result by Littelwood (see \cite{Li}
appendix p\ 295)

\begin{proposition}
\label{prop_in_A2n}Consider $\lambda $ and $\mu $ in $\mathcal{P}_{n}.$ Then:

\begin{enumerate}
\item  $[V_{n}^{B}(\lambda ):V_{2n}^{A}(\nu )]=\sum_{\gamma \in \mathcal{P}%
_{n}}c_{\gamma ,\lambda }^{\nu },$

\item  $[V_{n}^{C}(\lambda ):V_{2n}^{A}(\nu )]=\sum_{\gamma \in \mathcal{P}%
_{n}^{(1,1)}}c_{\gamma ,\lambda }^{\nu },$

\item  $[V_{n}^{D}(\lambda ):V_{2n}^{A}(\nu )]=\sum_{\gamma \in \mathcal{P}%
_{n}^{(2)}}c_{\gamma ,\lambda }^{\nu }.$
\end{enumerate}
\end{proposition}

\noindent The proposition below follows immediately from Theorem $A_{1}$ of 
\cite{K}. \ 

\begin{proposition}
Consider $\nu \in \mathcal{P}_{n}$ and $\lambda ^{+},\lambda ^{-}$ two
partitions such that $(\lambda ^{+},\lambda ^{-})$ has length $n.$ Then:

\begin{enumerate}
\item  $[V_{n}^{A}(\lambda ^{+},\lambda ^{-}):V_{n}^{B}(\nu )]=\sum_{\gamma
,\delta \in \mathcal{P}_{n}}c_{\gamma ,\delta }^{\nu }c_{\lambda
^{+},\lambda ^{-}}^{\delta },$

\item  $[V_{n}^{A}(\lambda ^{+},\lambda ^{-}):V_{n}^{C}(\nu )]=\sum_{\gamma
,\delta \in \mathcal{P}_{n}^{(2)}}c_{\gamma ,\delta }^{\nu }c_{\lambda
^{+},\lambda ^{-}}^{\delta },$

\item  $[V_{n}^{A}(\lambda ^{+},\lambda ^{-}):V_{n}^{D}(\nu )]=\sum_{\gamma
,\delta \in \mathcal{P}_{n}^{(1,1)}}c_{\gamma ,\delta }^{\nu }c_{\lambda
^{+},\lambda ^{-}}^{\delta }.$
\end{enumerate}
\end{proposition}

\noindent When $(\lambda^{+},\lambda^{-})=\lambda$ is a partition (that is $%
\lambda^{-}=\emptyset),$ we obtain the following dualities:

\begin{corollary}
\label{cor_mult}Consider $\lambda ,\nu $ two partitions of length $n,$ then

\begin{enumerate}
\item  $[V_{n}^{A}(\lambda ):V_{n}^{B}(\nu )]=[V_{n}^{B}(\lambda
):V_{2n}^{A}(\nu )]=\sum_{\gamma \in \mathcal{P}_{n}}c_{\gamma ,\lambda
}^{\nu },$

\item  $[V_{n}^{A}(\lambda ):V_{n}^{C}(\nu )]=[V_{n}^{D}(\lambda
):V_{2n}^{A}(\nu )]=\sum_{\gamma \in \mathcal{P}_{n}^{(2)}}c_{\gamma
,\lambda }^{\nu },$

\item  $[V_{n}^{A}(\lambda ):V_{n}^{D}(\nu )]=[V_{n}^{C}(\lambda
):V_{2n}^{A}(\nu )]=\sum_{\gamma \in \mathcal{P}_{n}^{(1,1)}}c_{\gamma
,\lambda }^{\nu }.$
\end{enumerate}
\end{corollary}

\subsection{ Kostka-Foulkes polynomials}

For any $w\in W_{B_{n}},$ the dot action of $w$ on $\beta\in\mathbb{Z}^{n}$
is defined by 
\begin{equation*}
w\circ\beta=w(\beta+\rho_{B_{n}})-\rho_{B_{n}}.
\end{equation*}
The $q$-analogue $\mathcal{P}_{q}^{B_{n}}$ of the Kostant partition function
corresponding to the root system $B_{n}$ is defined by the equality 
\begin{equation*}
\prod_{\alpha\in R_{B_{n}}^{+}}\dfrac{1}{1-qx^{\alpha}}=\sum_{\beta \in%
\mathbb{Z}^{n}}\mathcal{P}_{q}^{B_{n}}(\beta)x^{\beta}.
\end{equation*}
Note that $\mathcal{P}_{q}^{B_{n}}(\beta)=0$ if $\beta$ is not a linear
combination of positive roots of $R_{B_{n}}^{+}$with nonnegative
coefficients.\ We write similarly $\mathcal{P}_{q}^{C_{n}}$ and $\mathcal{P}%
_{q}^{D_{n}}$ for the $q$-partition functions associated respectively to the
root systems $C_{n}$ and $D_{n}$. Given $\lambda$ and $\mu$ two partitions
of length $n,$ the Kostka-Foulkes polynomials of types $B_{n},C_{n}$ and $%
D_{n}$ are then respectively defined by 
\begin{align*}
K_{\lambda,\mu}^{B_{n}}(q) & =\sum_{\sigma\in W_{B_{n}}}(-1)^{l(\sigma )}%
\mathcal{P}_{q}^{B_{n}}(\sigma(\lambda+\rho_{B_{n}})-(\mu+\rho_{B_{n}})), \\
K_{\lambda,\mu}^{C_{n}}(q) & =\sum_{\sigma\in W_{C_{n}}}(-1)^{l(\sigma )}%
\mathcal{P}_{q}^{C_{n}}(\sigma(\lambda+\rho_{C_{n}})-(\mu+\rho_{C_{n}})), \\
K_{\lambda,\mu}^{D_{n}}(q) & =\sum_{\sigma\in W_{D_{n}}}(-1)^{l(\sigma )}%
\mathcal{P}_{q}^{D_{n}}(\sigma(\lambda+\rho_{D_{n}})-(\mu+\rho_{D_{n}})).
\end{align*}
Set 
\begin{align*}
\widetilde{K}_{\lambda,\mu}^{B_{n}}(q) & =\sum_{\sigma\in\mathcal{S}%
_{n}}(-1)^{l(\sigma)}\mathcal{P}_{q}^{B_{n}}(\sigma(\lambda+\rho_{B_{n}})-(%
\mu+\rho_{n})), \\
\widetilde{K}_{\lambda,\mu}^{C_{n}}(q) & =\sum_{\sigma\in\mathcal{S}%
_{n}}(-1)^{l(\sigma)}\mathcal{P}_{q}^{C_{n}}(\sigma(\lambda+\rho_{C_{n}})-(%
\mu+\rho_{n})), \\
\widetilde{K}_{\lambda,\mu}^{D_{n}}(q) & =\sum_{\sigma\in\mathcal{S}%
_{n}}(-1)^{l(\sigma)}\mathcal{P}_{q}^{D_{n}}(\sigma(\lambda+\rho_{D_{n}})-(%
\mu+\rho_{n}))
\end{align*}
where $\rho_{n}=(n,...,1).\;$In \cite{lec}, we have proved that the
polynomials $\widetilde{K}_{\lambda,\mu}(q)$ are also Kostka-Foulkes
polynomials. More precisely we have:

\begin{lemma}
\label{lem_ktilde}Consider $\lambda ,\mu $ two partitions of length $n$ such
that $\left| \lambda \right| \geq \left| \mu \right| .\;$Let $k$ be any
integer such that $k\geq \frac{\left| \lambda \right| -\left| \mu \right| }{2%
}$.\ Then we have 
\begin{equation*}
\widetilde{K}_{\lambda ,\mu }(q)=K_{\lambda +k\kappa _{n},\mu +k\kappa
_{n}}(q)
\end{equation*}
where $\kappa _{n}=(1,...,1)\in \mathbb{Z}^{n}.$
\end{lemma}

\noindent\textbf{Remark: }Since $\sigma(\kappa_{n})=\kappa_{n}$ for any $%
\sigma\in\mathcal{S}_{n},$ we have $\widetilde{K}_{\lambda+k\kappa_{n},\mu+k%
\kappa_{n}}(q)=\widetilde{K}_{\lambda,\mu}(q)$ for any integer $k\geq0.$ So
we can extend the above definition of $\widetilde{K}_{\lambda,\mu}(q)$ for $%
\lambda$ and $\mu$ decreasing sequences of integers (positive or not).

\subsection{ The $q$-multiplicities $u_{\protect\lambda,\protect\mu}(q)$ and 
$U_{\protect\lambda,\protect\mu }(q)\label{sub_sec_def_U}$}

Set 
\begin{gather*}
\prod_{1\leq i<j\leq n}\dfrac{1}{1-q\frac{x_{i}}{x_{j}}}\prod_{1\leq r<s\leq
n}\dfrac{1}{1-\frac{q}{x_{i}x_{j}}}=\sum_{\beta\in\mathbb{Z}%
^{n}}f_{q}(\beta)x^{\beta}\text{ and } \\
\prod_{1\leq i<j\leq n}\dfrac{1}{1-q\frac{x_{i}}{x_{j}}}\prod_{1\leq r\leq
s\leq n}\dfrac{1}{1-\frac{q}{x_{i}x_{j}}}=\sum_{\beta\in\mathbb{Z}%
^{n}}F_{q}(\beta)x^{\beta}.
\end{gather*}
Given $\lambda$ and $\mu$ two partitions of length $n,$ let $u_{\lambda,\mu
}(q)$ and $U_{\lambda,\mu}(q)$ be the two polynomials defined by

\begin{equation*}
u_{\lambda,\mu}(q)=\sum_{\sigma\in\mathcal{S}_{n}}(-1)^{l(\sigma)}f_{q}(%
\sigma(\lambda+\rho_{n})-\mu-\rho_{n})\text{ and }U_{\lambda,\mu}(q)=\sum_{%
\sigma\in\mathcal{S}_{n}}(-1)^{l(\sigma)}F_{q}(\sigma(\lambda
+\rho_{n})-\mu-\rho_{n})
\end{equation*}
where $\rho_{n}=(n,...,1).$ Then $u_{\lambda,\mu}(q)$ and $%
U_{\lambda,\mu}(q) $ can be regarded as quantizations of tensor product
multiplicities \cite{lec}. Consider the tensor products 
\begin{gather*}
V_{\mu}^{B}=V^{B}(\mu_{1}\Lambda_{1})\otimes\cdot\cdot\cdot\otimes
V^{B}(\mu_{n}\Lambda_{1}),\text{ }V_{\mu}^{C}=V^{C}(\mu_{1}\Lambda_{1})%
\otimes \cdot\cdot\cdot\otimes V^{C}(\mu_{n}\Lambda_{1})\text{,} \\
V_{\mu}^{C}=V^{D}(\mu_{1}\Lambda_{1})\otimes\cdot\cdot\cdot\otimes
V^{D}(\mu_{n}\Lambda_{1})
\end{gather*}
and 
\begin{gather*}
W_{\mu}^{B}=W^{B}(\mu_{1}\Lambda_{1})\otimes\cdot\cdot\cdot\otimes
W^{B}(\mu_{n}\Lambda_{1}),\text{ }W_{\mu}^{C}=W^{C}(\mu_{1}\Lambda_{1})%
\otimes \cdot\cdot\cdot\otimes W^{C}(\mu_{n}\Lambda_{1})\text{,} \\
W_{\mu}^{D}=W^{D}(\mu_{1}\Lambda_{1})\otimes\cdot\cdot\cdot\otimes
W^{D}(\mu_{n}\Lambda_{1})
\end{gather*}
where for any $k\in\mathbb{N}$, $W(k_{1})=V(k\Lambda_{1})\oplus
V((k-2)\Lambda _{1})\oplus\cdot\cdot\cdot\oplus V((k\func{mod}%
2)\Lambda_{1}). $ Then we have the following proposition:

\begin{proposition}
\cite{lec}\label{prop_multi} Let $\lambda $ and $\mu $ be two partitions of
length $n.\;$Then

\begin{enumerate}
\item  $u_{\lambda ,\mu }(q)$ is a $q$-analogue of the multiplicity of the
representation $V(\lambda )$ in $V_{\mu },$

\item  $U_{\lambda ,\mu }(q)$ is a $q$-analogue of the multiplicity of the
representation $V(\lambda )$ in $W_{\mu }.$
\end{enumerate}
\end{proposition}

\noindent\textbf{Remarks:}

\noindent$\mathrm{(i):}$ It follows from the definition of $f_{q}$ and $%
F_{q} $ that $u_{\lambda,\mu}(q)=U_{\lambda,\mu}(q)=0$ if $\left|
\lambda\right| >\left| \mu\right| .$

\noindent$\mathrm{(ii):}$ When $q=1,$ we recover that the multiplicities of $%
V^{B}(\lambda),V^{C}(\lambda)$ and $V^{D}(\lambda)$ respectively in $V_{\mu
}^{B},$ $V_{\mu}^{C},$ and $V_{\mu}^{D}$ are equal \cite{K}.

\noindent$\mathrm{(iii):}$ In \cite{lec}, we have also obtained that $%
u_{\lambda,\mu}(q)$ and $U_{\lambda,\mu}(q)$ can be regarded as $q$%
-multiplicities in tensor product of column shaped representations.

\noindent$\mathrm{(iv):}$ Like the definition of $\widetilde{K}_{\lambda,\mu
}(q)$, the definitions of $u_{\lambda,\mu}(q)$ and $U_{\lambda,\mu}(q)$ can
also be extended for $\lambda$ and $\mu$ decreasing sequences of integers.

\noindent$\mathrm{(v):}$ Consider $\lambda,\mu\in\mathcal{P}_{n}$ and set $%
\lambda^{\#}=(\lambda_{1},...,\lambda_{n},0),$ $\mu^{\#}=(\mu_{1},...,\mu
_{n},0).\;$Then $u_{\lambda^{\#},\mu^{\#}}(q)=u_{\lambda,\mu}(q)$ and $%
U_{\lambda^{\#},\mu^{\#}}(q)=U_{\lambda,\mu}(q).$

\begin{theorem}
\cite{lec}\label{th_dual1} Consider $\lambda ,\mu $ two partitions of length 
$n$ and set $m=\max (\lambda _{1},\mu _{1})$.\ Then $\widehat{\lambda }%
=(m-\lambda _{n},...,m-\lambda _{1})$ and $\widehat{\mu }=(m-\mu
_{n},...,m-\mu _{1})$ are partitions of length $n$ and 
\begin{equation*}
u_{\lambda ,\mu }(q)=\widetilde{K}_{\widehat{\lambda },\widehat{\mu }%
}^{D_{n}}(q),\text{ }U_{\lambda ,\mu }(q)=\widetilde{K}_{\widehat{\lambda },%
\widehat{\mu }}^{C_{n}}(q).
\end{equation*}
\end{theorem}

\noindent Write $I$ for the involution defined on $\mathbb{Z}^{n}$ by $%
I(\beta_{1},...,\beta_{n})=(-\beta_{n},...,-\beta_{1}).$ For any decreasing
sequence of integers $\gamma,$ $I(\gamma)$ is also a decreasing sequence of
integers. By Theorem \ref{th_dual1}, this means that the correspondences 
\begin{equation}
u_{\lambda,\mu}(q)\longleftrightarrow\widetilde{K}_{I(\lambda),I(%
\mu)}^{D_{n}}(q)\text{ and }U_{\lambda,\mu}(q)\longleftrightarrow\widetilde{K%
}_{I(\lambda),I(\mu)}^{C_{n}}(q)  \label{corres}
\end{equation}
where $\lambda,\mu$ are decreasing sequence of integers can be interpreted
as a duality result for the $q$-multiplicities associated to the classical
root systems.

\subsection{Crystals of type $C_{n}$ and the one dimension sum $X_{\protect%
\lambda ,\protect\mu}(q)$ $\label{sub-sec-crts}$}

We have seen that $U_{\lambda,\mu}(q)$ can be regarded as a $q$-analogue of
the multiplicity of $V(\lambda)$ in$\frak{\;}W_{\mu}^{C}.\;$In \cite{Ok},
Hatayama, Kuniba, Okado and Takagi have introduced another quantification $%
X_{\lambda,\mu}(q)$ of this multiplicity based on the determination of the
combinatorial $R$-matrix of certain $U_{q}^{\prime}(C_{n}^{(1)})$-crystals $%
B_{k}^{C}$.\ Considered as the crystal graph of a $U_{q}(C_{n})$-module, $%
B_{k}^{C}$ is isomorphic to 
\begin{equation}
B^{C}(k\Lambda_{1})\oplus B^{C}((k-2)\Lambda_{1})\oplus\cdot\cdot\cdot\oplus
B^{C}(k\func{mod}2\Lambda_{1})  \label{sum}
\end{equation}
where for any $i\in\{k,k-2,...,k\func{mod}2\}$, $B^{C}(k\Lambda_{1})$ is the
crystal graph of the irreducible finite dimensional $U_{q}(C_{n})$-module of
highest weight $k\Lambda_{1}.\;$In \cite{KN}, Kashiwara and Nakashima have
obtained a natural labelling of the vertices of $B^{C}(k\Lambda_{1})$ by
one-row tableaux of length $k$ filled by letters of the alphabet 
\begin{equation*}
\mathcal{C}_{n}=\{1<\cdot\cdot\cdot<n-1<n<\overline{n}<\overline{n-1}%
<\cdot\cdot\cdot<\overline{1}\}
\end{equation*}
such that the letters increase from left to right (that is by semistandard
one-row tableaux on $\mathcal{C}_{n}$). Then the vertices of $B_{k}^{C}$ can
be depicted by one-row tableaux of length $k$ by adding $p$ pairs $(0,%
\overline{0})\ $to the tableaux appearing in the crystals $%
B^{C}((k-2p)\Lambda_{1})$ of the decomposition (\ref{sum}). Then by setting 
\begin{equation*}
\mathcal{C}_{n+1}=\{0<1<\cdot\cdot\cdot<n-1<n<\overline{n}<\overline
{n-1}<\cdot\cdot\cdot<\overline{1}<\overline{0}\}
\end{equation*}
the crystal $B_{k}^{C}$ can be regarded as a subcrystal of the $C_{n+1}$%
-crystal associated to the dominant weight $k\Lambda_{1}$ (labelled by the
one-row semistandard tableaux on $\mathcal{C}_{n+1}$ with length $k).$

\noindent Recall that the combinatorial $R$-matrix associated to crystals $%
B_{k}^{C}$ is equivalent to the description of the crystal graph
isomorphisms 
\begin{equation*}
\left\{ 
\begin{array}{c}
B_{l}^{C}\otimes B_{k}^{C}\overset{\simeq }{\rightarrow }B_{k}^{C}\otimes
B_{l}^{C} \\ 
b_{1}\otimes b_{2}\longmapsto b_{2}^{\prime }\otimes b_{1}^{\prime }
\end{array}
\right.
\end{equation*}
together with the energy function $H^{C}$ on $B_{l}^{C}\otimes B_{k}^{C}.\;$%
As proved in \cite{Ok}, this can done by using the insertion algorithm for $%
C_{n+1}$-tableaux of \cite{Ba} or \cite{lec2}. In the sequel we only need
the description of the energy function $H^{C}.$ Consider $b_{1}\in B_{l}^{C}$
and $b_{2}\in B_{k}^{C}$ and denote by $z=\min (\#0$ in $b_{1},$ $\#0$ in $%
b_{2})=\min (\#\overline{0}$ in $b_{1},$ $\#\overline{0}$ in $b_{2})$. Let $%
b_{1}^{\ast }$ and $b_{2}^{\ast }$ be the row tableaux obtained by erasing $%
z $ pairs of letters $(0,\overline{0})$ in $b_{1}$ and $b_{2}.$ Write $%
l^{\ast }$ and $k^{\ast }$ for the lengths of $b_{1}^{\ast }$ and $%
b_{2}^{\ast }.$ Denote by $P_{n+1}^{C}(b_{1}^{\ast }\otimes b_{2}^{\ast })$
the tableau obtained by inserting the row $b_{2}^{\ast }$ into the row $%
b_{1}^{\ast }$ following the $C_{n+1}$-insertion algorithm. Then $%
P_{n+1}^{C}(b_{1}^{\ast }\otimes b_{2}^{\ast })$ is a two-row $C_{n+1}$%
-tableau which contains $k^{\ast }+l^{\ast }$ letters. Since the plactic
relations for the root system $C_{n+1}$ are not homogeneous, pairs of
letters $(0,\overline{0})$ can appear in $P_{n+1}^{C}(b_{1}^{\ast }\otimes
b_{2}^{\ast }).$ Write $m$ for the length of the shortest row of $%
P_{n+1}^{C}(b_{1}^{\ast }\otimes b_{2}^{\ast }).$

\begin{proposition}
\cite{Ok} For any $b_{1}\otimes b_{2}$ in $B_{l}^{C}\otimes B_{k}^{C}$ we
have 
\begin{equation*}
H^{C}(b_{1}\otimes b_{2})=\min (l^{\ast },k^{\ast })-m.
\end{equation*}
\end{proposition}

\noindent The multiplicity of $V(\lambda)$ in $W_{\mu}^{C}$ is then equal to
the number of highest weight vertices of weight $\lambda$ in the crystal $%
B_{\mu}^{C}=B_{\mu_{1}}^{C}\otimes\cdot\cdot\cdot\otimes B_{\mu_{n}}^{C}$.\
Then the one dimension sum $X_{\lambda,\mu}(q)$ for $C_{n}^{(1)}$-crystals
is defined in \cite{Ok} by 
\begin{equation*}
X_{\lambda,\mu}(q)=\sum_{b\in E_{\lambda}}q^{H^{C}(b)}\text{ with }%
H^{C}(b)=\sum_{0\leq i<j\leq n}H^{C}(b_{i}\otimes b_{j}^{(i+1)})
\end{equation*}
where $E_{\lambda}$ is the set of highest weight vertices $b=b_{1}\otimes
\cdot\cdot\cdot\otimes b_{n}$ in $B_{\mu}^{C}$ of highest weight $\lambda,$ $%
b_{j}^{(i)}$ is determined by the crystal isomorphism 
\begin{equation}
\begin{array}{c}
B_{\mu_{i}}^{C}\otimes B_{\mu_{i+1}}^{C}\otimes B_{\mu_{i+2}}^{C}\otimes
\cdot\cdot\cdot\otimes B_{\mu_{j}}^{C}\rightarrow B_{\mu_{i}}^{C}\otimes
B_{\mu_{j}}^{C}\otimes B_{\mu_{i+1}}^{C}\cdot\cdot\cdot\otimes
B_{\mu_{j-1}}^{C} \\ 
b_{i}\otimes b_{i+1}\otimes\cdot\cdot\cdot\otimes b_{j}\rightarrow
b_{j}^{(i)}\otimes b_{i}^{\prime}\otimes\cdot\cdot\cdot\otimes
b_{j-1}^{\prime}
\end{array}
\label{def_isom}
\end{equation}
and for any $j=1,...,n,$ $H(b_{0}\otimes b_{j}^{(1)})$ is equal to the
number of letters $0$ in $b_{j}^{(1)}$.

\noindent When $\left| \lambda \right| =\left| \mu \right| $ the vertices of 
$E_{\lambda }$ contain only unbarred letters. In this case $%
H^{A}(b_{i}\otimes b_{j}^{(i+1)})=H^{C}(b_{i}\otimes b_{j}^{(i+1)})$ that is
the energy function of type $A$ defined on the vertices $b$ which do not
contain any barred letter is the restriction of that of type $C.\;$In \cite
{NY} Nakayashiki and Yamada have proved the equality $H^{A}(b)=\mathrm{ch}%
(Q(b))$ where $Q(b)$ is the semistandard tableau associated to $b$ by
generalizing the Robinson-Schensted correspondence and $\mathrm{ch}$ is the
charge defined by Lascoux and Sch\"{u}tzenberger in \cite{LS} and \cite{LSc1}%
.\ Recall that the charge statistic verifies 
\begin{equation}
K_{\lambda ,\mu }^{A_{n-1}}(q)=\sum_{T\in SST_{\mu }(\lambda )}q^{\mathrm{ch}%
(T)}  \label{LST}
\end{equation}
where $SST_{\mu }(\lambda )$ is the set of semistandard tableaux of shape $%
\lambda $ and weight $\mu .$ This implies that $X_{\lambda ,\mu
}(q)=K_{\lambda ,\mu }^{A_{n-1}}(q)$ when $\left| \lambda \right| =\left|
\mu \right| .$

\noindent Many computations suggest the following identities:

\begin{conjecture}
\label{conj}Consider $\lambda \in \mathcal{P}_{n}$.\ Then we have 
\begin{align*}
\mathrm{(i)}& :u_{\lambda ,(1^{n})}(q)=q^{\tfrac{\left| \mu \right| -\left|
\lambda \right| }{2}}X_{\lambda ,(1^{n})}(q), \\
\mathrm{(ii)}& :U_{\lambda ,\mu }(q)=q^{\left| \mu \right| -\left| \lambda
\right| }X_{\lambda ,\mu }(q)\text{ for any }\mu \in \mathcal{P}_{n}.
\end{align*}
\end{conjecture}

\noindent\textbf{Remarks:}

\noindent $\mathrm{(i)}$: A $q$-analogue for the multiplicity of $%
V_{n}^{C}(\lambda )$ in $W^{C}(\mu )$ can also be defined from rigged
configurations. The $X=M$ conjecture gives a simple relation called
fermionic formula between $X$ and $M.$ In \cite{OSS}, Okado, Schilling and
Shimozono have proved this conjecture when all the parts of $\mu $ are equal
to $1.$ In \ref{sec_proof}, we will prove $\mathrm{(i)}$ of Conjecture \ref
{conj} and $\mathrm{(ii)}$ when $\mu =(1,...,1)$. Thus by combining our
results with those of \cite{OSS} we obtain very simple relations between the
three different $q$-analogues for the multiplicity of $V_{n}^{C}(\lambda )$
in $W^{C}(1^{n})=V(1^{n})$.

\noindent$\mathrm{(ii)}$: By Theorem \ref{th_dual1}, the polynomials $%
U_{\lambda,\mu}(q)$ are Kostka-Foulkes polynomials.\ Thus $\mathrm{(ii)}$ of
Conjecture \ref{conj} implies that, up to a simple renormalization, the one
dimension sums $X_{\lambda,\mu}(q)$ are affine Kazhdan-Lusztig polynomials.

\section{Identities for the $q$-multiplicities}

\subsection{Decomposition of the $\widetilde{K}_{\protect\lambda,\protect\mu%
}(q)$ in terms of the $K_{\protect\lambda,\protect\mu}^{A_{n-1}}(q)$}

Denote by $\mathcal{P}_{q}^{A}$ the $q$-Kostant partition function defined
by 
\begin{equation*}
\prod_{1\leq i<j\leq n}\left( 1-q\frac{x_{i}}{x_{j}}\right) ^{-1}=\sum
_{\eta\in\mathbb{Z}^{n}}\mathcal{P}_{q}^{A}(\eta)x^{\eta}
\end{equation*}
Then for any $\lambda,\mu\in\mathcal{P}_{n},$ the Kostka-Foulkes polynomial $%
K_{\lambda,\mu}^{A_{n-1}}(q)$ is such that 
\begin{equation*}
K_{\lambda,\mu}^{A_{n-1}}(q)=\sum_{\sigma\in\mathcal{S}_{n}}(-1)^{l(\sigma )}%
\mathcal{P}_{q}^{A}(\sigma\circ\lambda-\mu).
\end{equation*}
Since $K_{\lambda+k\kappa_{n},\mu+k\kappa_{n}}^{A_{n-1}}(q)=K_{\lambda,\mu
}^{A_{n-1}}(q)$ for any positive integer $k$, the definition of $K_{\lambda
,\mu}^{A_{n-1}}(q)$ can be extended for $\lambda$ and $\mu$ decreasing
sequence of integers.\ By definition of the $q$-partition function $\mathcal{%
P}_{q}^{C}$ we must have 
\begin{equation*}
\sum_{\beta\in\mathbb{Z}^{n}}\mathcal{P}_{q}^{C}(\beta)x^{\beta}=\prod_{1%
\leq r\leq s\leq n}\left( 1-qx_{r}x_{s}\right) ^{-1}\prod_{1\leq i<j\leq
n}\left( 1-q\frac{x_{i}}{x_{j}}\right) ^{-1}.
\end{equation*}
Now we can write by (\ref{def_b}) 
\begin{equation}
\prod_{1\leq r\leq s\leq n}\left( 1-qx_{r}x_{s}\right) ^{-1}=\sum_{\delta\in
L_{C}}q^{\left| \delta\right| /2}c(\delta)x^{\delta}  \label{equal_fund}
\end{equation}
since the number of roots appearing in a decomposition of $\delta\in L_{C}$
as a sum of positive roots $\varepsilon_{r}+\varepsilon_{s}$ with $1\leq
r\leq s\leq n$ is always equal to $\left| \delta\right| /2.$ Thus we obtain 
\begin{equation}
\sum_{\beta\in\mathbb{Z}^{n}}\mathcal{P}_{q}^{C}(\beta)x^{\beta}=\sum_{\eta
\in\mathbb{Z}^{n}}\sum_{\delta\in L_{C}}q^{\left| \delta\right| /2}\mathcal{P%
}_{q}^{A}(\eta)c(\delta)x^{\delta+\eta}.  \label{P_C=P_A}
\end{equation}
By using similar arguments for $\mathcal{P}_{q}^{D}$ we derive the following
lemma:

\begin{lemma}
\label{lem_util}For any $\beta \in \mathbb{Z}^{n}$ we have 
\begin{equation*}
\mathcal{P}_{q}^{C}(\beta )=\sum_{\delta \in L_{C}^{\left| \beta \right|
}}q^{\left| \delta \right| /2}c(\delta )\mathcal{P}_{q}^{A}(\beta -\delta )%
\text{ and }\mathcal{P}_{q}^{C}(\beta )=\sum_{\delta \in L_{D}^{\left| \beta
\right| }}q^{\left| \delta \right| /2}d(\delta )\mathcal{P}_{q}^{A}(\beta
-\delta )
\end{equation*}
where $L_{C}^{\left| \beta \right| }=\{\delta \in L_{C},\left| \delta
\right| =\left| \beta \right| \}$ and $L_{D}^{\left| \beta \right|
}=\{\delta \in L_{D},\left| \delta \right| =\left| \beta \right| \}.$
\end{lemma}

\begin{proof}
\TEXTsymbol{>}From (\ref{P_C=P_A}) we derive the equality $\mathcal{P}%
_{q}^{C}(\beta )=\sum_{\eta +\delta =\beta }c(\delta )q^{\left| \delta
\right| /2}\mathcal{P}_{q}^{A_{n}}(\eta ).\;$Since $\mathcal{P}%
_{q}^{A_{n}}(\eta )=0$ when $\left| \eta \right| \neq 0,$ we can suppose $%
\left| \eta \right| =0$ and $\left| \delta \right| =\left| \beta \right| $
in the previous sum.\ Then $\delta \in L_{C}^{\left| \beta \right| }$ and
the result follows immediately.\ The proof for $\mathcal{P}_{q}^{D}(\beta )$
is similar.
\end{proof}

\noindent\textbf{Remark: }A similar result for the $q$-partition function $%
\mathcal{P}_{q}^{B}$ does not exit. Indeed the number of roots appearing in
a decomposition of $\delta\in L_{B}$ as a sum of positive roots $\varepsilon
_{r}+\varepsilon_{s}$ with $1\leq r<s\leq n$ and $\varepsilon_{i}$ with $%
1\leq i\leq n$ does not depend only of $\left| \delta\right| $ since $\left|
\varepsilon_{r}+\varepsilon_{s}\right| \neq\left| \varepsilon_{i}\right| .$

\begin{proposition}
\label{prop_dec_K_c}Consider $\lambda ,\mu \in \mathcal{P}_{n}$ such that $%
\left| \lambda \right| \geq \left| \mu \right| .$ Then we have:

\begin{enumerate}
\item  $\widetilde{K}_{\lambda ,\mu }^{C_{n}}(q)=q^{\tfrac{\left| \lambda
\right| -\left| \mu \right| }{2}}\sum_{\gamma \in \widetilde{\mathcal{P}}%
_{n}}\sum_{\sigma \in \mathcal{S}_{n}}(-1)^{l(\sigma )}c(\sigma \circ
\lambda -\gamma )K_{\gamma ,\mu }^{A_{n-1}}(q)$

\item  $\widetilde{K}_{\lambda ,\mu }^{D_{n}}(q)=q^{\tfrac{\left| \lambda
\right| -\left| \mu \right| }{2}}\sum_{\gamma \in \widetilde{\mathcal{P}}%
_{n}}\sum_{\sigma \in \mathcal{S}_{n}}(-1)^{l(\sigma )}d(\sigma \circ
\lambda -\gamma )K_{\gamma ,\mu }^{A_{n-1}}(q)$
\end{enumerate}

\noindent where $\widetilde{\mathcal{P}}_{n}=\{\gamma =(\gamma
_{1},...,\gamma _{n})\in \mathbb{Z}^{n},\gamma _{1}\geq \gamma _{2}\geq
\cdot \cdot \cdot \geq \gamma _{n}\}$.
\end{proposition}

\begin{proof}
We have 
\begin{equation*}
\widetilde{K}_{\lambda ,\mu }^{C_{n}}(q)=\sum_{\sigma \in \mathcal{S}%
_{n}}(-1)^{l(\sigma )}\mathcal{P}_{q}^{C_{n}}(\sigma (\lambda +\rho
_{n})-(\mu +\rho _{n})).
\end{equation*}
Hence from the previous lemma we derive 
\begin{equation*}
\widetilde{K}_{\lambda ,\mu }^{C_{n}}(q)=\sum_{\delta \in L_{C}^{\left|
\beta \right| }}c(\delta )q^{\left| \delta \right| /2}\sum_{\sigma \in 
\mathcal{S}_{n}}(-1)^{l(\sigma )}\mathcal{P}_{q}^{A_{n-1}}(\sigma (\lambda
+\rho _{n})-(\mu +\delta +\rho _{n}))
\end{equation*}
and 
\begin{equation*}
\widetilde{K}_{\lambda ,\mu }^{C_{n}}(q)=\sum_{\sigma \in \mathcal{S}%
_{n}}(-1)^{l(\sigma )}\sum_{\delta \in L_{C}^{\left| \beta \right|
}}c(\delta )q^{\left| \delta \right| /2}\mathcal{P}_{q}^{A_{n-1}}(\sigma
(\lambda +\rho _{n}-\sigma ^{-1}(\delta ))-(\mu +\rho _{n}))
\end{equation*}
For any $\sigma \in \mathcal{S}_{n},$ we have $\sigma ^{-1}(L_{C}^{\left|
\beta \right| })=L_{C}^{\left| \beta \right| }$ and $c(\delta )=c(\sigma
(\delta )).$ Thus we obtain 
\begin{equation}
\widetilde{K}_{\lambda ,\mu }^{C_{n}}(q)=\sum_{\sigma \in \mathcal{S}%
_{n}}(-1)^{l(\sigma )}\sum_{\delta \in L_{C}^{\left| \beta \right|
}}c(\delta )q^{\left| \delta \right| /2}\mathcal{P}_{q}^{A_{n-1}}(\sigma
(\lambda +\rho _{n}-\delta )-(\mu +\rho _{n}))=\sum_{\delta \in
L_{C}^{\left| \beta \right| }}c(\delta )q^{\frac{\left| \mu \right| -\left|
\lambda \right| }{2}}K_{\lambda -\delta ,\mu }^{A_{n-1}}(q).  \label{equK}
\end{equation}
\ Now $K_{\lambda -\delta ,\mu }^{A_{n-1}}(q)=0$ or there exits $\sigma \in 
\mathcal{S}_{n}$ and $\gamma \in \widetilde{\mathcal{P}}_{n}$ such that $%
\gamma =\sigma ^{-1}\circ (\lambda -\delta ).$ It follows that 
\begin{equation*}
\widetilde{K}_{\lambda ,\mu }^{C_{n}}(q)=\sum_{\sigma \in \mathcal{S}%
_{n}}(-1)^{l(\sigma )}\sum_{\gamma \in \widetilde{\mathcal{P}}%
_{n}}c((\lambda +\rho )-\sigma (\gamma +\rho ))K_{\gamma ,\mu }^{A_{n-1}}(q).
\end{equation*}
Since $c(\delta )=c(\sigma (\delta ))$ for any $\sigma \in \mathcal{S}_{n},$
and $\delta \in L_{C}$ we obtain the desired equality 
\begin{equation*}
\widetilde{K}_{\lambda ,\mu }^{C_{n}}(q)=\sum_{\gamma \in \widetilde{%
\mathcal{P}}_{n}}\sum_{\sigma \in \mathcal{S}_{n}}(-1)^{l(\sigma )}c(\sigma
(\lambda +\rho )-(\gamma +\rho ))K_{\gamma ,\mu }^{A_{n-1}}(q)
\end{equation*}
The proof is similar for $\widetilde{K}_{\lambda ,\mu }^{D_{n}}(q).$
\end{proof}

\noindent By Lemma \ref{lem_ktilde}, $\widetilde{K}_{\lambda,\mu}^{C_{n}}(q)$
and $\widetilde{K}_{\lambda,\mu}^{D_{n}}(q)$ are Kazhdan-Lusztig
polynomials.\ Thus they have nonnegative integer coefficients. This property
can also be obtained from the proposition below:

\begin{corollary}
\label{cor-dec_Ktilde}Consider $\lambda ,\mu \in \mathcal{P}_{n}.\;$For $k$
a sufficiently large integer we have:

\begin{enumerate}
\item  $\widetilde{K}_{\lambda ,\mu }^{C_{n}}(q)=q^{\tfrac{\left| \lambda
\right| -\left| \mu \right| }{2}}\sum_{\gamma \in \widetilde{\mathcal{P}}%
_{n}}[V_{n}^{A}(\gamma +k\kappa _{n}),V_{n}^{C}(\lambda +k\kappa
_{n})]K_{\gamma ,\mu }^{A_{n-1}}(q),$

\item  $\widetilde{K}_{\lambda ,\mu }^{D_{n}}(q)=q^{\tfrac{\left| \lambda
\right| -\left| \mu \right| }{2}}\sum_{\gamma \in \widetilde{\mathcal{P}}%
_{n}}[V_{n}^{A}(\gamma +k\kappa _{n}),V_{n}^{D}(\lambda +k\kappa
_{n})]K_{\gamma ,\mu }^{A_{n-1}}(q).$
\end{enumerate}
\end{corollary}

\begin{proof}
Consider $\gamma \in \widetilde{\mathcal{P}}_{n}$ and write $\gamma =(\gamma
^{+},\gamma ^{-})$ as in (\ref{jamma+-}).$\;$We have by Proposition \ref
{prop_mul_sum} $[V_{n}^{A}(\gamma ):V_{n}^{C}(\lambda )]=\sum_{w\in
W_{C_{n}}}(-1)^{l(w)}c(w\circ \lambda -\gamma ).$ Now if $k$ is sufficiently
large $\gamma +k\kappa _{n}$ is a partition and $[V_{n}^{A}(\gamma
):V_{n}^{C}(\lambda )]=\sum_{w\in W_{C_{n}}}(-1)^{l(w)}c(w\circ \lambda
-\gamma +w(k\kappa _{n})-k\kappa _{n}).$ Suppose that $w\notin \mathcal{S}%
_{n},$ then $\left| w(k\kappa _{n})-k\kappa _{n}\right| \leq -2k.$ Thus we
can choose $k$ such that $\gamma +k\kappa _{n}$ is a partition and $w\circ
\lambda -\gamma +w(k\kappa _{n})-k\kappa _{n}\notin L_{C}$ for any $w\in
W_{C_{n}}$ such that $w\notin \mathcal{S}_{n}.$ For such an integer $k$ we
have $c(w\circ \lambda -\gamma +w(k\kappa _{n})-k\kappa _{n})=0$.\ Then we
derive $1$ from $1$ of Proposition \ref{prop_dec_K_c}$.$ We prove $2$ by
using similar arguments.
\end{proof}

\bigskip

To obtain a decomposition of the polynomial $\widetilde{K}_{\lambda,\mu
}^{B_{n}}(q)$ as a sum of polynomials $K_{\gamma,\mu}^{A_{n-1}}(q),$ we need
first to obtain its decomposition in terms of the polynomials $\widetilde
{K}_{\nu,\mu}^{D_{n}}(q).$

\begin{proposition}
Consider $\lambda ,\mu \in \mathcal{P}_{n}$ such that $\left| \lambda
\right| \geq \left| \mu \right| $.\ Then for $k$ a sufficiently large
integer we have: 
\begin{equation*}
\widetilde{K}_{\lambda ,\mu }^{B_{n}}(q)=\sum_{\nu \in \widetilde{\mathcal{P}%
}_{n}}q^{\left| \lambda \right| -\left| \nu \right| }[V_{n}^{D}(\nu +k\kappa
_{n}),V_{n}^{B}(\lambda +k\kappa _{n})]\widetilde{K}_{\nu ,\mu }^{D_{n}}(q).
\end{equation*}
\end{proposition}

\begin{proof}
By definition of the partition functions $\mathcal{P}_{q}^{B_{n}}$ and $%
\mathcal{P}_{q}^{D_{n}}$ we can write 
\begin{equation*}
\sum_{\beta \in \mathbb{Z}^{n}}\mathcal{P}_{q}^{B_{n}}(\beta )x^{\beta
}=\sum_{\gamma \in \mathbb{Z}^{n}}\mathcal{P}_{q}^{D_{n}}(\gamma )x^{\gamma
}\prod_{k=1}^{n}\dfrac{1}{1-qx_{k}}.
\end{equation*}
We have $\prod_{k=1}^{n}\dfrac{1}{1-qx_{k}}=\sum_{\delta \in \mathbb{N}%
^{n}}q^{\left| \delta \right| }x^{\delta }.$ Thus we derive 
\begin{equation*}
\sum_{\beta \in \mathbb{Z}^{n}}\mathcal{P}_{q}^{B_{n}}(\beta )x^{\beta
}=\sum_{\gamma \in \mathbb{Z}^{n}}\sum_{\delta \in \mathbb{N}^{n}}q^{\left|
\delta \right| }\mathcal{P}_{q}^{D_{n}}(\gamma )x^{\gamma +\delta }.
\end{equation*}
This implies that $\mathcal{P}_{q}^{B_{n}}(\beta )=\sum_{\delta \in \mathbb{N%
}^{n}}q^{\left| \delta \right| }\mathcal{P}_{q}^{D_{n}}(\beta -\delta )$ and 
\begin{equation*}
\widetilde{K}_{\lambda ,\mu }^{B_{n}}(q)=\sum_{\sigma \in \mathcal{S}%
_{n}}(-1)^{l(\sigma )}\sum_{\delta \in \mathbb{N}^{n}}q^{\left| \delta
\right| }\mathcal{P}_{q}^{D_{n}}(\sigma (\lambda +\rho _{n})-(\mu +\delta
+\rho _{n})).
\end{equation*}
By using similar arguments to that of the proof of Proposition \ref
{prop_dec_K_c}, we obtain 
\begin{equation*}
\widetilde{K}_{\lambda ,\mu }^{B_{n}}(q)=q^{\left| \lambda \right| -\left|
\mu \right| }\sum_{\nu \in \widetilde{\mathcal{P}}_{n}}\sum_{\sigma \in 
\mathcal{S}_{n}}(-1)^{l(\sigma )}\mathbf{1}_{\mathbb{N}}(\sigma (\lambda
+\rho _{n})-(\nu +\rho _{n}))K_{\nu ,\mu }^{D_{n}}(q).
\end{equation*}
There exists a sufficiently large integer $k$ for which $w(\lambda +\rho
_{n}+k\kappa _{n})-(\nu +\rho _{n}+k\kappa _{n})\notin \mathbb{N}^{n}$ for
any $w\in W_{B_{n}}$ such that $w\notin \mathcal{S}_{n}$.\ Then it follows
from Lemma \ref{lem_dec_b_D} that $[V_{n}^{D}(\nu +k\kappa
_{n}):V_{n}^{B}(\lambda +k\kappa _{n})]=\sum_{\sigma \in \mathcal{S}%
_{n}}(-1)^{l(\sigma )}\mathbf{1}_{\mathbb{N}}(\sigma \circ \lambda -\nu )$
and the proposition is proved.
\end{proof}

\bigskip

\noindent\textbf{Remark:} From $2$ of Corollary \ref{cor-dec_Ktilde} and the
above proposition we obtain 
\begin{equation*}
\widetilde{K}_{\lambda,\mu}^{B_{n}}(q)=\sum_{\nu\in\widetilde{\mathcal{P}}%
_{n}}\sum_{\gamma\in\widetilde{\mathcal{P}}_{n}}q^{\frac{2\left|
\lambda\right| -\left| \nu\right| -\left| \mu\right| }{2}}[V_{n}^{D}(\nu+k%
\kappa_{n}),V_{n}^{B}(\lambda+k\kappa_{n})][V_{n}^{A}(\gamma+k\kappa
_{n}),V_{n}^{D}(\nu+k\kappa_{n})]K_{\gamma,\mu}^{A_{n-1}}(q)
\end{equation*}
which implies in particular that the polynomials $\widetilde{K}_{\lambda,\mu
}^{B_{n}}(q)$ have non negative coefficients.

\subsection{Decomposition of $u_{\protect\lambda,\protect\mu}(q)$ and $U_{%
\protect\lambda,\protect\mu}(q)$ in terms of the $K_{\protect\lambda,\protect%
\mu}^{A_{n-1}}(q)$}

By using similar arguments to the proof of Lemma \ref{lem_util}, we can
establish that for any $\beta\in\mathbb{Z}^{n}$ we have 
\begin{equation*}
f_{q}(\beta)=\sum_{\delta\in L_{C}^{\left| \beta\right| }}q^{\left|
\delta\right| /2}d(\delta)\mathcal{P}_{q}^{A}(\beta+\delta)\text{ and }%
F_{q}(\beta)=\sum_{\delta\in L_{D}^{\left| \beta\right| }}q^{\left|
\delta\right| /2}c(\delta)\mathcal{P}_{q}^{A}(\beta+\delta)
\end{equation*}
which implies the following proposition:

\begin{proposition}
\label{prp_dec_u}Consider $\lambda ,\mu \in \mathcal{P}_{n}$ such that $%
\left| \mu \right| \leq \left| \lambda \right| .\;$Then we have

\begin{enumerate}
\item  $U_{\lambda ,\mu }(q)=q^{\tfrac{\left| \mu \right| -\left| \lambda
\right| }{2}}\sum_{\nu \in \mathcal{P}_{n}}[V_{n}^{A}(\lambda
),V_{n}^{C}(\nu )]K_{\nu ,\mu }^{A_{n-1}}(q),$

\item  $u_{\lambda ,\mu }(q)=q^{\tfrac{\left| \mu \right| -\left| \lambda
\right| }{2}}\sum_{\nu \in \mathcal{P}_{n}}[V_{n}^{A}(\lambda
),V_{n}^{D}(\nu )]K_{\nu ,\mu }^{A_{n-1}}(q).$
\end{enumerate}
\end{proposition}

\begin{proof}
We proceed as in the proof of Proposition \ref{prop_dec_K_c} and (\ref{equK}%
) is replaced by 
\begin{equation*}
U_{\lambda ,\mu }(q)=\sum_{\delta \in L_{C}^{\left| \beta \right| }}c(\delta
)q^{\frac{\left| \mu \right| -\left| \lambda \right| }{2}}K_{\lambda +\delta
,\mu }^{A_{n-1}}(q).
\end{equation*}
This time $K_{\lambda +\delta ,\mu }^{A_{n-1}}(q)=0$ or there exits $\sigma
\in \mathcal{S}_{n}$ and $\nu \in \mathcal{P}_{n}$ such that $\nu =\sigma
^{-1}\circ (\lambda +\delta ).$ Note that $\nu $ can not have negative
coordinates for $\lambda +\delta $ have non negative coordinates. We deduce: 
\begin{equation*}
U_{\lambda ,\mu }(q)=q^{\frac{\left| \mu \right| -\left| \lambda \right| }{2}%
}\sum_{\sigma \in \mathcal{S}_{n}}(-1)^{l(\sigma )}\sum_{\nu \in \mathcal{P}%
_{n}}c(\sigma (\nu +\rho _{n})-(\lambda +\rho _{n}))K_{\gamma ,\mu
}^{A_{n-1}}(q).
\end{equation*}
We have seen in the proof of Corollary \ref{cor-dec_Ktilde} that for $k$ a
sufficiently large integer we have 
\begin{equation*}
\lbrack V_{n}^{A}(\lambda +k\kappa _{n}),V_{n}^{C}(\nu +k\kappa
_{n})]=\sum_{\sigma \in \mathcal{S}_{n}}(-1)^{l(\sigma )}c(\sigma \circ \nu
-\lambda ).
\end{equation*}
Since $\nu $ is a partition and $c_{\gamma ,\lambda +k\kappa _{n}}^{\nu
+k\kappa _{n}}=c_{\gamma ,\lambda }^{\nu }$ for any $k\in \mathbb{N}$, we
deduce from Corollary \ref{cor_mult} that $[V_{n}^{A}(\lambda
),V_{n}^{C}(\nu )]=[V_{n}^{A}(\lambda +k\kappa _{n}),V_{n}^{C}(\nu +k\kappa
_{n})]$ which proves $1.$ Assertion $2$ is obtained similarly.
\end{proof}

\subsection{Link with the polynomials $K_{\protect\lambda,\protect\mu%
}^{\diamondsuit}(q)$ of Shimozono and Zabrocki}

We deduce from Corollary \ref{cor_mult} and Proposition \ref{prp_dec_u}:

\begin{theorem}
\label{th-SZ}Consider $\lambda ,\mu \in \mathcal{P}_{n}$ such that $\left|
\mu \right| \geq \left| \lambda \right| \;$Then we have the following
equalities 
\begin{multline*}
\mathrm{(i)}:\text{ }u_{\lambda ,\mu }(q)=q^{\frac{\left| \mu \right|
-\left| \lambda \right| }{2}}\sum_{\nu \in \mathcal{P}_{n}}[V_{n}^{C}(%
\lambda ),V_{2n}^{A}(\nu )]K_{\nu ,\mu }^{A_{n-1}}(q)=q^{\frac{\left| \mu
\right| -\left| \lambda \right| }{2}}\sum_{\nu \in \mathcal{P}%
_{n}}[V_{n}^{A}(\lambda ),V_{n}^{D}(\nu )]K_{\nu ,\mu }^{A_{n-1}}(q) \\
=q^{\frac{\left| \mu \right| -\left| \lambda \right| }{2}}\sum_{\nu \in 
\mathcal{P}_{n}}\sum_{\gamma \in \mathcal{P}_{n}^{(1,1)}}c_{\gamma ,\lambda
}^{\nu }K_{\nu ,\mu }^{A_{n-1}}(q),
\end{multline*}
\begin{multline*}
\mathrm{(ii)}:\text{ }U_{\lambda ,\mu }(q)=q^{\frac{\left| \mu \right|
-\left| \lambda \right| }{2}}\sum_{\nu \in \mathcal{P}_{n}}[V_{n}^{D}(%
\lambda ),V_{2n}^{A}(\nu )]K_{\nu ,\mu }^{A_{n-1}}(q)=q^{\frac{\left| \mu
\right| -\left| \lambda \right| }{2}}\sum_{\nu \in \mathcal{P}%
_{n}}[V_{n}^{A}(\lambda ),V_{n}^{C}(\nu )]K_{\nu ,\mu }^{A_{n-1}}(q) \\
=q^{\frac{\left| \mu \right| -\left| \lambda \right| }{2}}\sum_{\nu \in 
\mathcal{P}_{n}}\sum_{\gamma \in \mathcal{P}_{n}^{(2)}}c_{\gamma ,\lambda
}^{\nu }K_{\nu ,\mu }^{A_{n-1}}(q).
\end{multline*}
\end{theorem}

\noindent By comparing the leftmost equalities of the above theorem with
equality (7.6) of \cite{SZ} we obtain 
\begin{equation}
K_{\lambda,\mu}^{(1,1)}(q)=u_{\lambda,\mu}(q^{2})\text{ and }K_{\lambda,\mu
}^{(2)}(q)=U_{\lambda,\mu}(q^{2})  \label{corees}
\end{equation}
where $K_{\lambda,\mu}^{(1,1)}(q)$ and $K_{\lambda,\mu}^{(2)}(q)$ are
polynomials defined by Shimozono and Zabrocki by using creating operators on
formal series. In particular the polynomials $K_{\lambda,\mu}^{(1,1)}(q)$
and $K_{\lambda,\mu}^{(2)}(q)$ are Kazhdan-Lusztig polynomials specialized
at $q^{2}.$

\bigskip

\noindent \textbf{Remark:} In \cite{SZ} the authors have also defined
another polynomial denoted $K_{\lambda ,\mu }^{(1)}(q)$ verifying 
\begin{equation}
K_{\lambda ,\mu }^{(1)}(q)=q^{\left| \mu \right| -\left| \lambda \right|
}\sum_{\nu \in \mathcal{P}_{n}}\sum_{\gamma \in \mathcal{P}_{n}}c_{\gamma
,\lambda }^{\nu }K_{\nu ,\mu }^{A_{n-1}}(q^{2}).  \label{def_K1}
\end{equation}
From the duality (\ref{corres}), it is tempting to introduce the polynomial $%
V_{\lambda ,\mu }(q)=\widetilde{K}_{I(\lambda ),I(\mu )}(q).\;$A similar
result than Proposition \ref{prp_dec_u} can not exist for $V_{\lambda ,\mu
}(q)$ (see the remark following Lemma \ref{lem_util}).\ Nevertheless, by
using Corollary \ref{cor_mult}, one can establish that 
\begin{equation*}
V_{\lambda ,\mu }(1)=\sum_{\nu \in \mathcal{P}_{n}}[V_{n}^{B}(\lambda
),V_{2n}^{A}(\nu )]K_{\nu ,\mu }^{A_{n-1}}(1)=\sum_{\nu \in \mathcal{P}%
_{n}}\sum_{\gamma \in \mathcal{P}_{n}}c_{\gamma ,\lambda }^{\nu }K_{\nu ,\mu
}^{A_{n-1}}(1).
\end{equation*}
Thus $K_{\lambda ,\mu }^{(1)}(1)=V_{\lambda ,\mu }(1).\;$However we have $%
K_{\lambda ,\mu }^{(1)}(q)\neq V_{\lambda ,\mu }(q^{2})$ in general.\ For
example if we take $\lambda =(1,0,0)$ and $\mu =(1,1,1)$ we obtain $%
K_{\lambda ,\mu }^{(1)}(q)=q^{8}+2q^{6}+2q^{4}+q^{2}$ and $V_{\lambda ,\mu
}(q^{2})=q^{10}+q^{8}+2q^{6}+q^{4}+q^{2}.$

\bigskip

Consider $\nu$ in $\mathcal{P}_{n}.$ For any standard tableau $\tau$ of
shape $\nu$ and weight $(1^{n}),$ let $\nu^{\prime}$ be the standard tableau
of shape $\nu^{\prime}$ and weight $(1^{n})$ obtained by reflecting $\tau$
among the diagonal. Then one can verify that $\mathrm{ch}(\tau^{\prime})=%
\tfrac{n(n-1)}{2}-\mathrm{ch}(\tau)$ which by (\ref{LST}) implies the
identity: 
\begin{equation}
K_{\nu^{\prime},(1^{n})}^{A_{n-1}}(q)=q^{\frac{n(n-1)}{2}}K_{%
\nu,(1^{n})}^{A_{n-1}}(q^{-1}).  \label{dual(1)n}
\end{equation}
The following proposition will be useful in Section \ref{sec_proof}.

\begin{proposition}
\label{prop_conj}Consider $\lambda \in \mathcal{P}_{n}$ such that $n\geq
\left| \lambda \right| \;$Then we have 
\begin{equation*}
\mathrm{(i):}\text{ }U_{\lambda ^{\prime },(1^{n})}(q)=q^{\frac{n(n-1)}{2}%
+n-\left| \lambda \right| }u_{\lambda ,(1^{n})}(q^{-1})\text{ and }\mathrm{%
(ii):}\text{ }u_{\lambda ^{\prime },(1^{n})}(q)=q^{\frac{n(n-1)}{2}+n-\left|
\lambda \right| }U_{\lambda ,(1^{n})}(q^{-1}).
\end{equation*}
\end{proposition}

\begin{proof}
By $\mathrm{(ii)}$ of Theorem \ref{th-SZ} we have 
\begin{equation}
U_{\lambda ^{\prime },(1^{n})}(q)=q^{\tfrac{n-\left| \lambda \right| }{2}%
}\sum_{\nu \in \mathcal{P}_{n}}\sum_{\gamma \in \mathcal{P}%
_{n}^{(2)}}c_{\gamma ,\lambda ^{\prime }}^{\nu }K_{\nu
,(1^{n})}^{A_{n-1}}(q).  \label{aux}
\end{equation}
Since $K_{\nu ,(1^{n})}^{A_{n-1}}(q)=0$ when $\left| \nu \right| \neq n,$ we
can suppose that $\nu $ belongs to $\widehat{\mathcal{P}}_{n}=\{\nu \in 
\mathcal{P}_{n},$ $\left| \nu \right| =n\}$ in the above sum$.\;$For any $%
\nu $ in $\widehat{\mathcal{P}}_{n},$ we have $c_{\gamma ,\lambda }^{\nu }=0$
unless $\left| \gamma \right| =n-\left| \lambda \right| .\;$So we can
suppose that $\gamma $ belongs to $\widehat{\mathcal{P}}_{n}^{(2)}=\{\gamma
\in \mathcal{P}_{n}^{(2)},$ $\left| \gamma \right| \leq n\}$ in (\ref{aux})$%
.\;$The map $\Gamma :\nu \longmapsto \nu ^{\prime }$ is an involution of $%
\widehat{\mathcal{P}}_{n}.\;$Moreover $\Gamma (\widehat{\mathcal{P}}%
_{n}^{(2)})=\widehat{\mathcal{P}}_{n}^{(1,1)}=\{\gamma \in \mathcal{P}%
_{n}^{(1,1)},$ $\left| \gamma \right| \leq n\}.$ Thus we can write 
\begin{equation*}
U_{\lambda ^{\prime },(1^{n})}(q)=q^{\tfrac{n-\left| \lambda \right| }{2}%
}\sum_{\nu \in \widehat{\mathcal{P}}_{n}}\sum_{\gamma \in \widehat{\mathcal{P%
}}_{n}^{(1,1)}}c_{\gamma ^{\prime },\lambda ^{\prime }}^{\nu ^{\prime
}}K_{\nu ^{\prime },(1^{n})}^{A_{n-1}}(q).
\end{equation*}
Now since $c_{\gamma ^{\prime },\lambda ^{\prime }}^{\nu ^{\prime
}}=c_{\gamma ,\lambda }^{\nu },$ we derive by (\ref{dual(1)n}) 
\begin{equation*}
U_{\lambda ^{\prime },(1^{n})}(q)=q^{\frac{n(n-1)}{2}+\tfrac{n-\left|
\lambda \right| }{2}}\sum_{\nu \in \widehat{\mathcal{P}}_{n}}\sum_{\gamma
\in \widehat{\mathcal{P}}_{n}^{(1,1)}}c_{\gamma ,\lambda }^{\nu }K_{\nu
,(1^{n})}^{A_{n-1}}(q^{-1}).
\end{equation*}
Finally the equality $U_{\lambda ^{\prime },(1^{n})}(q)=q^{\frac{n(n-1)}{2}%
+n-\left| \lambda \right| }u_{\lambda ,(1^{n})}(q^{-1})$ is deduced from $%
\mathrm{(i)}$ of Theorem \ref{th-SZ}.\ We obtain $\mathrm{(ii)}$ similarly.
\end{proof}

\bigskip

\noindent\textbf{Remark: }By introducing for the root system $A_{n-1}$
generalized Kostka-Foulkes polynomials $K_{\nu,R}(q)$ where $%
R=(R_{1},...,R_{n})$ is a sequence of rectangular partitions, Schilling and
Warnarr \cite{SW} have proved the equality $K_{\nu^{\prime},R^{%
\prime}}(q)=q^{\left\| R\right\| }K_{\nu,R^{\prime}}(q^{-1})$ where $%
R^{\prime}=(R_{1}^{\prime },...,R_{l}^{\prime})$ and $\left\| R\right\|
=\sum_{i<j}\left| R_{i}\cap R_{j}\right| $ which can be considered as a
generalization of (\ref{dual(1)n}). In \cite{SZ}, Shimozono and Zabrocki
have also defined their polynomials $K_{\lambda,R}^{\diamondsuit}(q)$ when $%
R $ is a sequence of rectangular partitions. By (\ref{corees}), the above
proposition can also be regarded as a Corollary of Proposition 28 of \cite
{SZ}.

\section{Proof of Conjecture \ref{conj} when $\protect\mu=(1,...,1)\label%
{sec_proof}$}

\subsection{The crystals $B_{\Xi}$}

For any\ integer $l\geq 0,$ let $B^{A}(l)$ be the crystal graph of the
irreducible finite dimensional $U_{q}(sl_{2n})$-module of highest weight $%
l\Lambda _{1}.\;$In the sequel we choose to label the vertices of $B^{A}(1)$
by the letters of $\mathcal{C}_{n},$ that is we identify $B^{A}(1)$ with the
crystal 
\begin{equation*}
1\overset{1}{\rightarrow }2\cdot \cdot \cdot \cdot \rightarrow n-1\overset{%
n-1}{\rightarrow }n\overset{n}{\rightarrow }\overline{n}\overset{n+1}{%
\rightarrow }\overline{n-1}\overset{n+2}{\rightarrow }\cdot \cdot \cdot
\cdot \rightarrow \overline{2}\overset{2n-1}{\rightarrow }\overline{1}.
\end{equation*}
Recall that the crystal graph $B^{C}(1)$ has been identified in 2.5 with 
\begin{equation*}
1\overset{1}{\rightarrow }2\cdot \cdot \cdot \cdot \rightarrow n-1\overset{%
n-1}{\rightarrow }n\overset{n}{\rightarrow }\overline{n}\overset{n-1}{%
\rightarrow }\overline{n-1}\overset{n-2}{\rightarrow }\cdot \cdot \cdot
\cdot \rightarrow \overline{2}\overset{1}{\rightarrow }\overline{1}.
\end{equation*}
Thus the crystals $B^{A}(1)$ and $B^{C}(1)$ have the same vertices. For any
partition $\mu \in \mathcal{P}_{n},$ set $B_{(\mu )}^{A}=B^{A}(\mu
_{1})\otimes \cdot \cdot \cdot \otimes B^{A}(\mu _{n})$ and $B_{(\mu
)}^{C}=B^{C}(\mu _{1}\Lambda _{1})\otimes \cdot \cdot \cdot \otimes
B^{C}(\mu _{n}\Lambda _{1})$ (note that $B_{(\mu )}^{C}\neq B_{\mu }^{C}$
defined in \ref{sub-sec-crts})$.\;$Then $B_{(\mu )}^{A}$ and $B_{(\mu )}^{C}$
have also the same vertices.\ Nevertheless their crystal structure are
distinct and their decompositions in connected components do not coincide.

\noindent Denote by $H^{A}$ the energy function associated to $%
B^{A}(l)\otimes B^{A}(k).$ Then for any $b_{1}\otimes b_{2}$ belonging to $%
B^{A}(l)\otimes B^{A}(k),$ we have $H^{A}(b_{1}\otimes b_{2})=\min (l,k)-m$
where $m$ is the length of the shortest row of the semistandard tableau $%
P^{A}(b_{1}\otimes b_{2})$ obtained by inserting the row $b_{2}$ in the row $%
b_{1}$ following the column bumping algorithm for semistandard tableaux.
Given $b=b_{1}\otimes \cdot \cdot \cdot \otimes b_{n}\in B_{(\mu )}^{A},$
set 
\begin{equation*}
H^{A}(b)=\sum_{1\leq i<j\leq n}H^{A}(b_{i}\otimes b_{j}^{(i+1)})
\end{equation*}
where the vertices $b_{j}^{(i+1)}$ are defined as in (\ref{def_isom}). In
the sequel we need the following result due to Nakayashiki and Yamada:

\begin{theorem}
\label{th_NY}\cite{NY} Consider $\nu $ and $\mu $ two partitions of $%
\mathcal{P}_{n}$.\ Then 
\begin{equation*}
K_{\nu ,\mu }^{A_{n-1}}(q)=\sum_{b\in G_{\nu }}q^{H^{A}(b)}
\end{equation*}
where $G_{\nu }$ is the set of highest weight vertices of weight $\nu $ in $%
B_{(\mu )}^{A}.$
\end{theorem}

\noindent \textbf{Remark:} It is possible to show that $P^{A}(b_{1}\otimes
b_{2})=P_{n+1}^{C}(b_{1}\otimes b_{2})$ (see \ref{sub-sec-crts}) for any $%
b_{1}\otimes b_{2}$ belonging to $B^{A}(l)\otimes B^{A}(k)$ if and only if $%
l=k=1.\;$Moreover if we choose $l\geq 2$ and $k\geq 2,$ these two tableaux
can have distinct shapes.$\;$For example, by taking $n=2,$ $l=k=2,$ $%
b_{1}=b_{2}=\bar{2}2,$ we obtain 
\begin{equation*}
P^{A}(b_{1}\otimes b_{2})= 
\begin{tabular}{|l|ll}
\hline
$\mathtt{2}$ & $\mathtt{2}$ & \multicolumn{1}{|l|}{$\mathtt{\bar{2}}$} \\ 
\hline
$\mathtt{\bar{2}}$ &  &  \\ \cline{1-1}
\end{tabular}
\text{ and }P_{n+1}^{C}(b_{1}\otimes b_{2})= 
\begin{tabular}{|l|l|}
\hline
$\mathtt{1}$ & $\mathtt{\bar{2}}$ \\ \hline
$\mathtt{2}$ & $\mathtt{\bar{1}}$ \\ \hline
\end{tabular}
.
\end{equation*}
Hence the two statistics $H^{A}$ and $H^{C}$ do not coincide in general on
the vertices of $B_{(\mu )}^{A}.$

\bigskip

\noindent Suppose now that $\mu =(1,...,1).\;$Then any vertex $b$ of $%
B_{(1^{n})}^{A}$ or $B_{(1^{n})}^{C}$ can be written 
\begin{equation*}
b=x_{1}\otimes \cdot \cdot \cdot \otimes x_{n}
\end{equation*}
where $x_{1},...,x_{n}$ are letters of $\mathcal{C}_{n}.\;$We have 
\begin{equation*}
H^{A}(b)=H^{C}(b)=\sum_{i=1}^{n-1}(n-i)H(x_{i}\otimes x_{i+1})
\end{equation*}
where $H(x_{i}\otimes x_{i+1})=1$ if $x_{i}\geq x_{i+1}$ and $H(x_{i}\otimes
x_{i+1})=0$ otherwise.

\noindent To each vertex $b,$ we associate the $(n-1)$-tuple $\Xi(b)=(\xi
_{1},...,\xi_{n-1})$ such that for any $i=1,...,n-1,$ $\xi_{i}=0$ if $%
x_{i}<x_{i+1}$ and $\xi_{i+1}=1$ otherwise.\ For any $(n-1)$-tuple $\Xi$
with coordinates equal to $0$ or $1,$ set 
\begin{equation*}
B_{\Xi}=\{b\in B_{(1^{n})}^{A},\Xi(b)=\Xi\}.
\end{equation*}
Then the statistics $H^{A}$ and $H^{C}$ are invariant on the vertices of $%
B_{\Xi}$ and we have 
\begin{equation*}
H^{A}(b)=H^{C}(b)=\theta(b)=\sum_{i=1}^{n-1}(n-i)\xi_{i}\text{ for any }b\in
B_{\Xi}\text{ with }\Xi=(\xi_{1},...,\xi_{n-1}).
\end{equation*}

\begin{lemma}
\label{lem_Bxhi}Let $\Xi =(\xi _{1},...,\xi _{n-1})$ be a $(n-1)$-tuple with
coordinates equal to $0$ or $1.\;$Then the set $B_{\Xi }$ has a structure of 
$A_{2n-1}$-crystal and a structure of $C_{n}$-crystal.
\end{lemma}

\begin{proof}
Consider $\widetilde{K}$ a Kashiwara crystal operator for $U_{q}(sl_{2n-1})$
or $U_{q}(sp_{2n})$ and $x,y$ two letters of $\mathcal{C}_{n}.\;$Set $%
\widetilde{K}(x\otimes y)=x^{\prime }\otimes y^{\prime }.\;$From the
description of the crystals $B^{A}(1)^{\otimes 2}$ and $B^{C}(1)^{\otimes 2}$
given by Kashiwara and Nakashima in \cite{KN} we derive the equivalence: 
\begin{equation*}
x\leq y\Longleftrightarrow x^{\prime }\leq y^{\prime }.
\end{equation*}
This implies that $B_{\Xi }$ is stable under the action of any Kashiwara
crystal operator.\ Thus $B_{\Xi }$ has a structure of $A_{2n-1}$-crystal and
a structure of $C_{n}$-crystal.
\end{proof}

\subsection{The $X=u$ conjecture when $\protect\mu=(1,...,1)$}

\begin{theorem}
\label{proof1}For any partition $\lambda $ of length $n$ we have 
\begin{equation*}
u_{\lambda ,(1^{n})}(q)=q^{\tfrac{n-\left| \lambda \right| }{2}}X_{\lambda
,(1^{n})}(q).
\end{equation*}
\end{theorem}

\begin{proof}
Denote by $\frak{E}$ the set of the $(n-1)$-tuples $\Xi $ with coordinates
equal to $0$ or $1$.\ By Lemma \ref{lem_Bxhi}, for any $\Xi \in \frak{E}$, $%
B_{\Xi }$ has a structure of $A_{2n-1}$-crystal and a structure of $C_{n}$%
-crystal.\ Let us denote respectively by $B_{\Xi }^{A}$ and $B_{\Xi }^{C}$
these two crystals. There exists a $U_{q}(gl_{2n})$-module $M_{\Xi }^{A}$
whose crystal is isomorphic to $B_{\Xi }^{A}.\;$Similarly there exists a $%
U_{q}(sp_{2n})$-module $M_{\Xi }^{C}$ whose crystal is isomorphic to $B_{\Xi
}^{A}.\;$Recall that the weight $\mathrm{wt}^{A}(b)$ of a vertex $b\in
B_{\Xi }^{A}$ is equal to $(d_{1},...,d_{n},d_{n+1},...,d_{2n})$ where for
any $i\in \{1,...,2n\},$ $d_{i}$ is the number of letters $i$ in $b.$
Similarly the weight $\mathrm{wt}^{C}(b)$ of a vertex $b\in B_{\Xi }^{C}$ is
equal to $(\delta _{1},...,\delta _{n})$ where for any $i\in \{1,...,2n\},$ $%
\delta _{i}$ is the number of letters $i$ in $b$ minus the number of letters 
$\overline{i}.$ Then the characters of the modules $M_{\Xi }^{A}$ and $%
M_{\Xi }^{C}$ verify 
\begin{equation*}
\mathrm{char}(M_{\Xi }^{A})(x_{1},...,x_{n})=\sum_{b\in B_{\Xi }^{A}}x^{%
\mathrm{wt}^{A}(b)}\text{ and }\mathrm{char}(M_{\Xi
}^{C})(x_{1},...,x_{n})=\sum_{b\in B_{\Xi }^{C}}x^{\mathrm{wt}^{C}(b)}.
\end{equation*}
Let $N_{\Xi }^{A}$ and $N_{\Xi }^{C}$ be two representations respectively of 
$GL_{2n}$ and $Sp_{2n}$ such that $\mathrm{char}(N_{\Xi }^{A})=\mathrm{char}%
(M_{\Xi }^{A})$ and $\mathrm{char}(N_{\Xi }^{C})=\mathrm{char}(M_{\Xi
}^{C}). $ Then, by definition of $\mathrm{wt}^{A}(b)$ and $\mathrm{wt}%
^{C}(b),$ we have 
\begin{equation*}
\mathrm{char}(N_{\Xi }^{C})(x_{1},...,x_{n})=\mathrm{char}(N_{\Xi
}^{A})(x_{1},...,x_{n},\dfrac{1}{x_{n}},...,\dfrac{1}{x_{1}})
\end{equation*}
that is $\mathrm{ch}(N_{\Xi }^{C})$ is obtained by changing $x_{n+i}$ in $%
\dfrac{1}{x_{i}}$ in $\mathrm{char}(N_{\Xi }^{A}).$ Thus $\mathrm{char}%
(N_{\Xi }^{C})=\mathrm{char}(N_{\Xi }^{A}\downarrow _{Sp_{2n}}^{GL_{2n}})$
where $N_{\Xi }^{A}\downarrow _{Sp_{2n}}^{GL_{2n}}$ is the restriction of $%
N_{\Xi }^{A}$ to the action of $Sp_{2n}.$ This implies that $N_{\Xi }^{C}$
and $N_{\Xi }^{A}\downarrow _{Sp_{2n}}^{GL_{2n}}$ are isomorphic
representations of $Sp_{2n}.$

\noindent Write $E_{\lambda }^{\Xi }=\{b\in B_{\Xi }^{C},$ $b\in E_{\lambda
}\}$ for the set of highest weight vertices of weight $\lambda $ in $B_{\Xi
}^{C}.$ Then we obtain: 
\begin{equation*}
\mathrm{card}(E_{\lambda }^{\Xi })=[V_{n}^{C}(\lambda ),N_{\Xi
}^{C}]=[V_{n}^{C}(\lambda ),N_{\Xi }^{A}]
\end{equation*}
where $[V_{n}^{C}(\lambda ),N_{\Xi }^{A}]$ is the multiplicity of $%
V_{n}^{C}(\lambda )$ in $N_{\Xi }^{A}\downarrow _{Sp_{2n}}^{GL_{2n}}.\;$%
Recall that $X_{\lambda ,(1^{n})}(q)=\sum_{b\in E_{\lambda }}q^{H^{C}(b)}$
where $E_{\lambda }$ is the set of highest weight vertices of weight $%
\lambda $ in $B_{(1^{n})}^{C}.$ By Lemma \ref{lem_Bxhi} and the equality
above we can write 
\begin{equation}
X_{\lambda ,(1^{n})}(q)=\sum_{\Xi \in \frak{E}}\sum_{b\in E_{\lambda }^{\Xi
}}q^{H^{C}(b)}=\sum_{\Xi \in \frak{E}}\mathrm{card}(E_{\lambda }^{\Xi
})q^{\theta _{\Xi }}=\sum_{\Xi \in \frak{E}}[V_{n}^{C}(\lambda ),N_{\Xi
}^{A}]q^{\theta _{\Xi }}  \label{eq1}
\end{equation}
where for any $\Xi =(\xi _{1},...,\xi _{n-1})$ in $\frak{E,}$ $\theta _{\Xi
}=\sum_{i=1}^{n-1}(n-i)\xi _{i}.$

\noindent By definition of the representations $N_{\Xi }^{A}$ we have 
\begin{equation}
\tbigoplus_{\Xi \in \frak{E}}N_{\Xi }^{A}\simeq V_{2n}^{A}(\Lambda
_{1})^{\otimes n}  \label{dec_A}
\end{equation}
as $GL_{2n}$-representations. The irreducible components of $B_{(1^{n})}^{A}$
are indexed by standard tableaux with letters in $\{1,...,n\}.$ For any $%
b\in B_{(1^{n})}^{A}$ denote by $Q(b)$ the recording tableau (which is a
standard tableau with letters in $\{1,...,n\})$ associated to $b$ by the
Robinson-Schensted correspondence. The $Q$-symbol yields a one to one
correspondence between the highest weight vertices of $B_{(1^{n})}^{A}$
(thus the irreducible components of $V_{2n}^{A}(\Lambda _{1})^{\otimes n})$
and the set $ST$ of standard tableaux with letters in $\{1,...,n\}.$ For any 
$\tau \in ST,$ denote by $V_{n}^{A}(\tau )$ the irreducible component of $%
V_{2n}^{A}(\Lambda _{1})^{\otimes n}$ associated to $\tau .$ Then $%
V_{2n}^{A}(\tau )\simeq V_{2n}^{A}(\nu )$ where $\nu $ is the partition
corresponding to the shape of $\tau .$ Write $ST_{\Xi }$ for the set of
standards tableaux $\tau \in ST$ such that there exists a highest weight
vertex $b\in B_{\Xi }^{A}$ (thus for the $A_{2n-1}$-structure of graph) with 
$Q(b)=\tau $.\ Then from (\ref{dec_A}) the sets $ST_{\Xi }$, $\Xi \in \frak{E%
}$ are disjoints and $\cup _{\frak{E}}ST_{\Xi }=$ $ST$.\ For any $\Xi \in 
\frak{E,}$ we have 
\begin{equation*}
N_{\Xi }^{A}=\tbigoplus_{\tau \in ST_{\Xi }}V_{2n}^{A}(\tau ).
\end{equation*}
Thus we derive from (\ref{eq1}) that 
\begin{equation*}
X_{\lambda ,(1^{n})}(q)=\sum_{\Xi \in \frak{E}}\sum_{\tau \in ST_{\Xi
}}[V_{n}^{C}(\lambda ),V_{2n}^{A}(\tau )]q^{\theta _{\Xi }}.
\end{equation*}
Now the multiplicity $[V_{n}^{C}(\lambda ),V_{2n}^{A}(\tau )]$ depends only
on the shape $\nu $ of the standard tableau $\tau .\;$This means that we
have the equality $[V_{n}^{C}(\lambda ),V_{2n}^{A}(\tau
)]=[V_{n}^{C}(\lambda ),V_{2n}^{A}(\nu )]$ for any $\tau $ of shape $\nu \in 
\mathcal{P}_{n}.$ We deduce that 
\begin{equation*}
X_{\lambda ,(1^{n})}(q)=\sum_{\Xi \in \frak{E}}\sum_{\nu \in \mathcal{P}%
_{n}}\sum_{\tau \in ST_{\Xi }^{\nu }}[V_{n}^{C}(\lambda ),V_{2n}^{A}(\nu
)]q^{\theta _{\Xi }}
\end{equation*}
where $ST_{\Xi }^{\nu }=\{\tau \in ST_{\Xi }$ of shape $\nu \}.$ So we
obtain 
\begin{equation*}
X_{\lambda ,(1^{n})}(q)=\sum_{\nu \in \mathcal{P}_{n}}[V_{n}^{C}(\lambda
),V_{2n}^{A}(\nu )]\sum_{\Xi \in \frak{E}}\mathrm{card}(ST_{\Xi }^{\nu
})q^{\theta _{\Xi }}.
\end{equation*}
But we have $\sum_{\Xi \in \frak{E}}\mathrm{card}(ST_{\Xi }^{\nu })q^{\theta
_{\Xi }}=K_{\nu ,(1^{n})}^{A_{n-1}}(q)$ since $\theta _{\Xi }=H^{A}(b)$ for
any highest weight vertex $b$ of weight $\nu $ in $B_{\Xi }^{A}.$ Thus 
\begin{equation*}
X_{\lambda ,(1^{n})}(q)=\sum_{\nu \in \mathcal{P}_{n}}[V_{n}^{C}(\lambda
),V_{2n}^{A}(\nu )]K_{\nu ,(1^{n})}^{A_{n-1}}(q).
\end{equation*}
Finally by Theorem \ref{th-SZ}, we obtain 
\begin{equation*}
q^{\frac{n-\left| \lambda \right| }{2}}X_{\lambda ,(1^{n})}(q)=u_{\lambda
,(1^{n})}(q)
\end{equation*}
and the Theorem is proved.
\end{proof}

\bigskip

\noindent \textbf{Remark: }One can define, from the crystal $B_{(\mu
)}^{C},\ $the sum 
\begin{equation*}
Y_{\lambda ,\mu }(q)=\sum_{b\in F_{\lambda }}q^{H^{C}(b)}
\end{equation*}
where $F_{\lambda }$ is the set of highest weight vertices $b=b_{1}\otimes
\cdot \cdot \cdot \otimes b_{n}$ in $B_{(\mu )}^{C}$ of weight $\lambda .$
However, the polynomial $Y_{\lambda ,\mu }(q)$ is not a one dimension sum
related to an affine root system and the identity 
\begin{equation*}
u_{\lambda ,\mu }(q)=q^{\tfrac{\left| \mu \right| -\left| \lambda \right| }{2%
}}Y_{\lambda ,\mu }(q)
\end{equation*}
is false in general. This is in particular the case for $\mu =(2,2,2)$ and $%
\lambda =(2,0,0).$ Nevertheless, from (\ref{corees}) and Conjecture 32 of 
\cite{SZ}, one can conjecture that the $q$-multiplicities $u_{\lambda ,\mu
}(q)$ are equal, up to the multiplication by a power of $q,$ to the one
dimension sums related to the affine root system $A_{2n-1}^{(2)}$.

\subsection{The $X=U$ conjecture when $\protect\mu =(1,...,1)$}

We want to establish the equality $U_{\lambda,(1^{n})}(q)=q^{n-\left|
\lambda\right| }X_{\lambda,(1^{n})}(q)$.\ By Proposition \ref{prop_conj} we
know that 
\begin{equation*}
U_{\lambda,(1^{n})}(q)=q^{\frac{n(n-1)}{2}+n-\left| \lambda\right|
}u_{\lambda^{\prime},(1^{n})}(q^{-1}).
\end{equation*}
Thus by Theorem \ref{proof1} we obtain 
\begin{equation*}
U_{\lambda,(1^{n})}(q)=q^{\frac{n(n-1)}{2}+\frac{n-\left| \lambda\right| }{2}%
}X_{\lambda^{\prime},(1^{n})}(q^{-1})
\end{equation*}
and it suffices to prove the equality 
\begin{equation*}
X_{\lambda^{\prime},(1^{n})}(q^{-1})=q^{-\tfrac{n(n-1)}{2}+\frac{n-\left|
\lambda\right| }{2}}X_{\lambda,(1^{n})}(q)
\end{equation*}
which is equivalent to 
\begin{equation}
X_{\lambda^{\prime},(1^{n})}(q)=q^{\tfrac{n(n-1)}{2}-\frac{n-\left|
\lambda\right| }{2}}X_{\lambda,(1^{n})}(q^{-1}).  \label{eq-X}
\end{equation}
In \cite{lec2} we have introduced a Robinson-Schensted type correspondence
for the vertices of $B_{(1^{n})}^{C}.\;$In particular we have obtained a one
to one correspondence between the highest weight vertices of $b\in
B_{(1^{n})}^{C}$ and the oscillating tableaux of length $n.\;$Recall that an
oscillating tableau of length $n$ is a sequence $Q=(Q_{1},...,Q_{n})$ of
Young diagrams such that $Q_{i}$ and $Q_{i+1}$ differs by exactly one box
(that is $Q_{k+1}/Q_{k}= 
\begin{tabular}{|l|}
\hline
\\ \hline
\end{tabular}
$\ or $Q_{k}/Q_{k+1}= 
\begin{tabular}{|l|}
\hline
\\ \hline
\end{tabular}
).\;$We denote by $Q(b)$ the oscillating tableau associated to the highest
weight vertex $b\in B_{(1^{n})}^{C}$ under this correspondence. More
precisely set $b=x_{1}\otimes\cdot\cdot\cdot\otimes x_{n}$.\ Then $Q(b)$ is
defined recursively as follows: 
\begin{equation}
Q_{1}= 
\begin{tabular}{|l|}
\hline
\\ \hline
\end{tabular}
\text{ and }Q_{i+1}=x_{i+1}\rightarrow Q_{i}  \label{def_Q}
\end{equation}
where $x_{i+1}\rightarrow Q_{i}$ is the Young diagram obtained from $Q_{i}$
by adding a box on the $k$-th row of $Q_{i}$ if $x_{i+1}=k\in\{1,...,n\}$
and by deleting a box on the $k$-th row of $Q_{i}$ if $x_{i+1}=\overline{k}%
\in\{\overline{1},...,\overline{n}\}.$ Given any highest weight vertex $b,$
it is easy to verify that $Q(b)$ is an oscillating tableau of length $n.$

\noindent Suppose $Q(b)=(Q_{1},...,Q_{n}).\;$Then $Q^{\prime}=(Q_{1}^{\prime
},...,Q_{n}^{\prime}),$ where for any $i\in\{1,...,n\}$ $Q_{i}^{\prime}$ is
the conjugate diagram of $Q_{i},$ is also an oscillating tableau \ There
exists a highest weight vertex $b^{\prime}\in B_{(1^{n})}^{C}$ such that $%
Q(b^{\prime})=Q^{\prime}.$ Moreover if the weight of $b$ is equal to $%
\lambda,$ then the weight of $b^{\prime}$ is equal to $\lambda^{\prime}.$ To
prove our conjecture we need the two following technical lemmas:

\begin{lemma}
\label{lem1}Suppose that $b=x_{1}\otimes \cdot \cdot \cdot \otimes x_{n}$ is
a highest weight vertex of $B_{(1^{n})}^{C}$ and write $b^{\prime
}=x_{1}^{\prime }\otimes \cdot \cdot \cdot \otimes x_{n}^{\prime }$.\ Then:

\begin{enumerate}
\item  for any $i=1,...,n,$ the two letters $x_{i}$ and $x_{i}^{\prime }$
are simultaneously barred or unbarred,

\item  $H(x_{i}^{\prime }\otimes x_{i+1}^{\prime })=\left\{ 
\begin{tabular}{l}
$1-H(x_{i}\otimes x_{i+1})$ if $x_{i}$ and $x_{i+1}$ are simultaneoulsly
barred or unbarred \\ 
$H(x_{i}\otimes x_{i+1})$ otherwise
\end{tabular}
\right. .$
\end{enumerate}
\end{lemma}

\begin{proof}
We obtain $1$ immediately from the definition (\ref{def_Q}) of $Q$ and $%
Q^{\prime }.$

\noindent As usual we enumerate the rows (resp.\ the columns) of the Young
diagrams from top to bottom (resp.\ from left to right).

\noindent Suppose that $x_{i}=p$ and $x_{i+1}=q$ with $p,q\in \{1,...,n\}$.\
Then $Q_{i+1}$ is obtained by adding first a box in the $p$-th row of $%
Q_{i-1}$ to give $Q_{i},$ next a box in the $q$-th row of $Q_{i}$.\ Thus $%
Q_{i+1}^{\prime }$ is obtained by adding first a box in the $p$-th column of 
$Q_{i-1}^{\prime },$ next a box in the $q$-th column of $Q_{i}^{\prime }$ if 
$p\neq q,$ by adding two boxes in the $p$-th column of $Q_{i-1}^{\prime }$
otherwise.\ This implies that we have $x_{i}^{\prime }\geq x_{i+1}^{\prime }$
when $p<q$ (i.e.\ $x_{i}<x_{i+1})$ and $x_{i}^{\prime }<x_{i+1}^{\prime }$
when $p\geq q$ (i.e.\ $x_{i}\geq x_{i+1}).$

\noindent Now suppose that $x_{i}=\overline{p}$ and $x_{i+1}=\overline{q}$
with $p,q\in \{1,...,n\}$.\ Then $Q_{i+1}$ is obtained by deleting first a
box in the $p$-th row of $Q_{i-1}$ to give $Q_{i},$ next a box in the $q$-th
row of $Q_{i}$.\ Thus $Q_{i+1}^{\prime }$ is obtained by deleting first a
box in the $p$-th column of $Q_{i-1}^{\prime },$ next a box in the $q$-th
column of $Q_{i}^{\prime }$ if $p\neq q,$ by deleting two boxes in the $p$%
-th column of $Q_{i-1}^{\prime }$ otherwise.\ This implies that we have $%
x_{i}^{\prime }\geq x_{i+1}^{\prime }$ when $\overline{p}<\overline{q}$
(i.e.\ $x_{i}<x_{i+1})$ and $x_{i}^{\prime }<x_{i+1}^{\prime }$ when $%
\overline{p}\geq \overline{q}$ (i.e.\ $x_{i}\geq x_{i+1}).$

\noindent So in all cases we obtain $H(x_{i}^{\prime }\otimes
x_{i+1}^{\prime })=1-H(x_{i}\otimes x_{i+1})$ when $x_{i}$ and $x_{i+1}$ are
simultaneously barred or unbarred.

\noindent The equality $H(x_{i}^{\prime }\otimes x_{i+1}^{\prime
})=H(x_{i}\otimes x_{i+1})$ when $x_{i}$ and $x_{i+1}$ are not
simultaneously barred or unbarred follows from $1$ since a barred letter is
always greater than an unbarred one.
\end{proof}

\bigskip

\noindent For any vertex $b=x_{1}\otimes\cdot\cdot\cdot\otimes x_{n}\in
B_{(1^{n})}^{C},$ set $Z_{b}=\{i\in\{1,...,n-1\},$ $x_{i}$ and $x_{i+1}$ are
not simultaneously barred or unbarred$\}$

\begin{lemma}
\label{lem2}With the above notation we have 
\begin{equation*}
\sum_{i\in Z_{b}}(n-i)(1-2H(x_{i}\otimes x_{i+1}))=\dfrac{n-\left| \lambda
\right| }{2}
\end{equation*}
for any highest weight vertex $b=x_{1}\otimes \cdot \cdot \cdot \otimes
x_{n}\in B_{(1^{n})}^{C}$ of weight $\lambda .$
\end{lemma}

\begin{proof}
Observe first that we have 
\begin{equation*}
1-2H(x_{i}\otimes x_{i+1})=\left\{ 
\begin{tabular}{l}
$1$ if $x_{i}\leq x_{i+1}$ \\ 
$-1$ otherwise
\end{tabular}
\right. .
\end{equation*}
Since $b$ is a highest weight vertex, it follows from the description of the
action of the Kashiwara operators \cite{Ka} on a tensor product of crystals
that for any $i\in \{1,...,n\},$ $b_{i}=x_{1}\otimes \cdot \cdot \cdot
\otimes x_{i}$ is a highest weight vertex.$\;$In particular we must have $%
b_{1}=x_{1}=1.$ Thus we obtain 
\begin{equation}
\sum_{i\in Z_{b}}(1-2H(x_{i}\otimes x_{i+1}))=\left\{ 
\begin{tabular}{l}
$0$ if $x_{n}$ is unbarred \\ 
$1$ otherwise
\end{tabular}
\right. .  \label{s_aux}
\end{equation}
To prove the lemma we proceed by induction on $n.$ When $n=1,$ we have $b=1,$
$Z_{b}=0$ and $\left| \lambda \right| =1.\;$Hence the lemma is true. Now
suppose the lemma proved for any highest weight vertex of $B_{(1^{n-1})}^{C}$
and consider $b=x_{1}\otimes \cdot \cdot \cdot \otimes x_{n}\in
B_{(1^{n})}^{C}$ a highest weight vertex of weight $\lambda .$ Set $%
s_{b}=\sum_{i\in Z_{b}}(n-i)(1-2H(x_{i}\otimes x_{i+1})).$ We have $%
Z_{b}=Z_{b_{n-1}}$ if $x_{n-1}$ and $x_{n}$ are simultaneously barred or
unbarred and $Z_{b}=Z_{b_{n-1}}\cup \{n-1\}$ otherwise.\ Write $\widetilde{%
\lambda }$ for the weight of the highest weight vertex $b_{n-1}=x_{1}\otimes
\cdot \cdot \cdot \otimes x_{n-1}$.\ In the two cases $Z_{b}=Z_{b_{n-1}}$
and $Z_{b}=Z_{b_{n-1}}\cup \{n-1\}$ we derive by using the induction
hypothesis 
\begin{equation*}
s_{b}=\frac{n-1-\left| \widetilde{\lambda }\right| }{2}+\sum_{i\in
Z_{b}}(1-2H(x_{i}\otimes x_{i+1})).
\end{equation*}
When $x_{n}$ is unbarred, we have $\left| \widetilde{\lambda }\right|
=\left| \lambda \right| -1.\;$Hence $\tfrac{n-1-\left| \widetilde{\lambda }%
\right| }{2}=\tfrac{n-\left| \lambda \right| }{2}$ and we obtain $s_{b}=%
\tfrac{n-\left| \lambda \right| }{2}$ by (\ref{s_aux}).$\;$When $x_{n}$ is
barred, we have $\left| \widetilde{\lambda }\right| =\left| \lambda \right|
+1$ thus $\tfrac{n-1-\left| \widetilde{\lambda }\right| }{2}=\tfrac{n-\left|
\lambda \right| }{2}-1$ and by (\ref{s_aux}) we also derive $s_{b}=\tfrac{%
n-\left| \lambda \right| }{2}$ which proves the lemma.
\end{proof}

\begin{theorem}
\label{th_U}For any partition $\lambda $ of length $n$ we have 
\begin{equation*}
U_{\lambda ,(1^{n})}(q)=q^{n-\left| \lambda \right| }X_{\lambda ,(1^{n})}(q).
\end{equation*}
\end{theorem}

\begin{proof}
Consider $b=x_{1}\otimes \cdot \cdot \cdot \otimes x_{n}$ a highest weight
vertex of $B_{(1^{n})}^{C}$ of weight $\lambda $ and set $b^{\prime
}=x_{1}^{\prime }\otimes \cdot \cdot \cdot \otimes x_{n}^{\prime }.$ Then 
\begin{equation*}
H(b^{\prime })=\sum_{i=1}^{n-1}(n-i)H(x_{i}^{\prime }\otimes x_{i+1}^{\prime
})=\sum_{i\in Z_{b}}(n-i)H(x_{i}\otimes x_{i+1})+\sum_{i\notin
Z_{b}}(n-i)(1-H(x_{i}\otimes x_{i+1}))
\end{equation*}
by Lemma \ref{lem1}. Thus we have 
\begin{equation*}
H(b^{\prime })=\sum_{i=1}^{n-1}(n-i)(1-H(x_{i}\otimes x_{i+1}))+\sum_{i\in
Z_{b}}(n-i)(2H(x_{i}\otimes x_{i+1})-1)=\frac{n(n-1)}{2}-\dfrac{n-\left|
\lambda \right| }{2}-H(b)
\end{equation*}
where the last equality follows from Lemma \ref{lem2}. Finally we derive 
\begin{equation*}
X_{\lambda ^{\prime },(1^{n})}(q)=\sum_{b\in F_{\lambda }}q^{H(b^{\prime
})}=q^{\frac{n(n-1)}{2}-\frac{n-\left| \lambda \right| }{2}}\sum_{b\in
F_{\lambda }}q^{-H(b^{\prime })}=q^{\frac{n(n-1)}{2}-\frac{n-\left| \lambda
\right| }{2}}X_{\lambda ,(1^{n})}(q^{-1})
\end{equation*}
which by (\ref{eq-X}) proves the Theorem.
\end{proof}

\bigskip

\noindent\textbf{Remark: }Theorem \ref{th_U} can be regarded as an analogue
for the affine root system $C_{n}^{(1)}$ of Theorem \ref{th_NY} when $%
\mu=(1,...,1).$

\subsection{Appendix on the one dimension sums}

In \cite{KKM}, Kang, Kashiwara and Misra have defined crystals $B_{l}^{\phi
} $ for quantum affine algebras associated to the affine root systems $\phi $
of types $A_{n}^{(1)},C_{n}^{(1)},A_{2n-1}^{(2)},A_{2n}^{(2)},$ $%
D_{n+1}^{(2)},$ $B_{n}^{(1)}$ and $D_{n}^{(1)}$.\ We give below the
decomposition of the crystals $B_{l}^{\phi }$ as classical crystals for each
affine root system $\phi $: 
\begin{equation}
\begin{tabular}{|c|c|c|c|c|c|}
\hline
$C_{n}^{(\overset{\text{ \ \ }}{1})}$ & $A_{2n-1}^{(2)}$ & $A_{2n}^{(2)}$ & $%
D_{n+1}^{(2)}$ & $B_{n}^{(1)}$ & $D_{n}^{(1)}$ \\ \hline
$\underset{k=0,k\equiv l\func{mod}2}{\overset{\overset{\text{ \ \ }}{l}}{%
\tbigoplus }}B^{C}(k\Lambda _{1})$ & $B^{C}(k\Lambda _{1})$ & $\underset{k=0%
}{\overset{l}{\tbigoplus }}B^{C}(k\Lambda _{1})$ & $\underset{k=0}{\overset{l%
}{\tbigoplus }}B^{B}(k\Lambda _{1})$ & $B^{B}(k\Lambda _{1})$ & $%
B^{D}(k\Lambda _{1})$ \\ \hline
\end{tabular}
\label{dec_class}
\end{equation}
where $B^{C}(k\Lambda _{1})$, $B^{D}(k\Lambda _{1})$ and $B^{D}(k\Lambda
_{1})$ are the crystal graphs of the finite dimensional modules of highest
weight $k\Lambda _{1}$ respectively for $U_{q}(\frak{sp}_{2n})$, $U_{q}(%
\frak{so}_{2n+1})$ and $U_{q}(\frak{so}_{2n}).$ Given any partition $\mu \in 
\mathcal{P}_{n},$ set $B_{\mu }^{\phi }=B_{\mu _{1}}^{\phi }\otimes \cdot
\cdot \cdot \otimes B_{\mu _{n}}^{\phi }.$ Then, by using the energy
functions explicited in \cite{Ok} and \cite{Ok2}, one can define a statistic 
$H^{\phi }$ on the vertices of $B_{\mu }^{\phi }$ from which it is possible
to calculate the one dimension sum $X_{\lambda ,\mu }^{\phi }(q).\;$In
particular, with our notation, we have 
\begin{equation*}
X_{\lambda ,\mu }^{C_{n}^{(1)}}(q)=X_{\lambda ,\mu }(q)
\end{equation*}
since the polynomial $X_{\lambda ,\mu }(q)$ is a one dimension sum for the
affine root system $C_{n}^{(1)}.\;$Denote by $m$ the number of nonzero parts
of $\mu .\;$One proves that the highest weight vertices of weight $\lambda $
in $B_{\mu }^{B_{n}^{(1)}}$ (considered as a $B_{n}$-crystal) are in one to
one correspondence with those of $B_{\mu }^{A_{2n-1}^{(2)}}$ (considered as
a $C_{n}$-crystal)$.\;$Similarly when $m<n$ the highest weight vertices of
weight $\lambda $ in $B_{\mu }^{D_{n}^{(1)}}$ (considered as a $D_{n}$%
-crystal) are also in one to one correspondence with those of $B_{\mu
}^{A_{2n-1}^{(2)}}.\;$Moreover the statistics $%
H^{A_{2n-1}^{(2)}},H^{B_{n}^{(1)}}$ and $H^{D_{n}^{(1)}}$ are conserved via
these correspondences. Thus we obtain by (\ref{dec_class}) the equality 
\begin{equation*}
X_{\lambda ,\mu }^{A_{2n}^{(2)}}(q)=X_{\lambda ,\mu }^{D_{n+1}^{(2)}}(q)
\end{equation*}
for any $\mu \in \mathcal{P}_{n}.\;$Similarly, when $m<n$, we have 
\begin{equation}
X_{\lambda ,\mu }^{A_{2n-1}^{(2)}}(q)=X_{\lambda ,\mu
}^{B_{n}^{(1)}}(q)=X_{\lambda ,\mu }^{D_{n}^{(1)}}(q).  \label{one_dim_sums}
\end{equation}
From Theorems \ref{th_dual1} and \ref{th_U} we obtain:

\begin{corollary}
\label{cor_final}For any $\lambda \in \mathcal{P}_{n}$ we have 
\begin{equation*}
X_{\lambda ,\mu }^{C_{n}^{(1)}}(q)=q^{-\frac{\left| \mu \right| -\left|
\lambda \right| }{2}}U_{\lambda ,\mu }(q)=q^{-\frac{\left| \mu \right|
-\left| \lambda \right| }{2}}\widetilde{K}_{\widehat{\lambda },\widehat{\mu }%
}^{C_{n}}(q)
\end{equation*}
when $\mu =(1,...,1)\in \mathcal{P}_{n}.$
\end{corollary}

\noindent When $\phi$ is one of the two affine root systems $A_{2n}^{(2)}$
or $D_{n+1}^{(2)},$ the structure of classical crystal of $B_{\mu}^{\phi}$
implies that the one dimension sum $X_{\lambda,\mu}^{\phi}(q)$ cannot be
naturally related to the $q$-multiplicities $u_{\lambda,\mu}(q)$ or $%
U_{\lambda,\mu}(q)$ (see \ref{sub_sec_def_U}). However they can be expressed
in terms of the polynomials $K_{\lambda,\mu}^{(1)}(q)$ \cite{sh}.

\bigskip

\noindent \textbf{Note: }\textit{While revising this work, the author have
been informed that, in a paper in preparation \cite{sh}, Shimozono obtains a
proof of the }$X=K$ \textit{conjecture for tensor product of the ``symmetric
power'' Kirilov-Reshetikin modules for nonexceptional affine algebras\ of
type }$D_{n+1}^{(2)},A_{2n}^{(2)},C_{n}^{(1)},A_{2n}^{(2)}$ or $%
A_{2n}^{(2)\dagger }$.\ \textit{With our convention for the definition of
the one dimension sums (which is that of \cite{Ok}) this result can be
reformulated by writing} 
\begin{equation*}
K_{\lambda ,\mu }^{(2)}(q)=q^{\left| \mu \right| -\left| \lambda \right|
}X_{\lambda ,\mu }^{C_{n}^{(1)}}(q^{2})=q^{\left| \mu \right| -\left|
\lambda \right| }X_{\lambda ,\mu }^{A_{2n}^{(2)\dagger }}(q^{2})\text{ and }%
K_{\lambda ,\mu }^{(1)}(q)=q^{\left| \mu \right| -\left| \lambda \right|
}X_{\lambda ,\mu }^{A_{2n}^{(2)}}(q^{2})=q^{\left| \mu \right| -\left|
\lambda \right| }X_{\lambda ,\mu }^{D_{n}^{(2)}}(q^{2})
\end{equation*}
\textit{for any }$\mu \in \mathcal{P}_{n}$\textit{.\ By (\ref{corees}) the
first equality above establishes }$\mathrm{(ii)}$ of \textit{Conjecture \ref
{conj} since }$X_{\lambda ,\mu }(q)$ \textit{is a one dimension sum for }$%
C_{n}^{(1)}$\textit{-crystals. Thus assertion }$2$\textit{\ of Corollary \ref
{cor_final} holds for any }$\mu \in \mathcal{P}_{n}.$


\bigskip

\begin{thebibliography}{99}
\bibitem{Ba}  \textsc{T. H.Baker }\textit{An} \textit{insertion scheme for }$%
C_{n}$\textit{\ crystals}, in M.\ Kashiwara and T. Miwa, eds., Physical
Combinatorics, Birkh\"{a}user, Boston, 2000, \textbf{191}: 1-48.

\bibitem{FH}  \textsc{W.\ Fulton, J. Harris, }\textit{Representation theory}%
, Graduate Texts in Mathematics, Springer-Verlag.

\bibitem{GW}  \textsc{G. Goodman, N. R Wallach, }\textit{Representation
theory and invariants of the classical groups}, Cambridge University Press.

\bibitem{Ok}  \textsc{G. Hatayama, A. Kuniba, M. Okado, T. Takagi,}\textit{\
Combinatorial }$R$ matrices for a family of crystals: $C_{n}^{(1)}$ and $%
A_{2n-1}^{(2)}$ cases, Physical Combinatorics edited by M Kashiwara and T
Miwa, Birkhauser, 105-139 (2000).

\bibitem{Ok2}  \textsc{H. Hatayama, A. Kuniba, M. Okado, T. Takagi,}\textit{%
\ Combinatorial }$R$ matrices for a family of crystals: $%
B_{n}^{(1)},D_{n}^{(1)},A_{2n}^{(2)}$ and $D_{n+1}^{(2)}$ cases, Journal of
Algebra, \textbf{247}, 577-615 (2002).

\bibitem{HKOTY}  \textsc{G. Hatayama, A. Kuniba, M. Okado, T. Takagi, Y.
Yamada,}\textit{\ Remarks on fermionic formula, }in N. Jing and K.\ C.\
Misra, eds.\ Recent Developments in Quantum Affine Algebras and Related
Topics, Contemporary Mathematics \textbf{248}, AMS, Providence, 243-291,
(1999).

\bibitem{KKM}  \textsc{S-J. Kang, M. Kashiwara, K-C. Misra, }\textit{Crystal
bases of Verma modules for quantum affine Lie algebras, }Compositio.\ Math.%
\textit{\ \textbf{92}, 299-325 (1994).}

\bibitem{Ka}  \textsc{M. Kashiwara}, \textit{On crystal bases of the }$q$%
\textit{-analogue of universal enveloping algebras}, Duke Math. J, \textbf{63%
} (1991), 465-516.

\bibitem{KN}  \textsc{M. Kashiwara, T. Nakashima,} \textit{Crystal graphs
for representations of the }$q$\textit{-analogue of classical Lie algebras},
Journal of Algebra, \textbf{165}, 295-345 (1994).

\bibitem{K}  \textsc{K. Koike,}\textit{\ }\textsc{I. Terada,}\textit{\ Young
diagrammatic methods for the representations theory of the classical groups
of type }$B_{n},C_{n}$ and $D_{n},$ Journal of Algebra, \textbf{107},
466-511 (1987).

\bibitem{KT}  \textsc{K. Koike,}\textit{\ }\textsc{I. Terada,}\textit{\
Restriction Rules for }$GL,SO,Sp,$ Adv.\ in Math., \textbf{79}, 104-135
(1990).

\bibitem{LS}  \textsc{A. Lascoux, M-P. Sch\"{u}tzenberger, }\ \textit{Le mono%
}$\mathit{\ddot{\imath}}$\textit{de plaxique}, in non commutative structures
in algebra and geometric combinatorics A. de Luca Ed., Quaderni della
Ricerca Scientifica del C.N.R., Roma, (1981)\textit{.}

\bibitem{LSc1}  \textsc{A. Lascoux, M-P. Sch\"{u}tzenberger, }\ \textit{Sur
une conjecture de H.O Foulkes}, CR Acad Sci Paris, \textbf{288}, 95-98
(1979).

\bibitem{lec}  \textsc{C. Lecouvey,} \textit{A duality between }$q$\textit{%
-multiplicities in tensor products and }$q$\textit{-multiplicities of
weights for the root systems }$B,C$\textit{\ or }$D$ (submitted), ArXiv:
CO/04007522.

\bibitem{lec2}  \textsc{C. Lecouvey,} \textit{Schensted-type correspondence,
Plactic Monoid and Jeu de Taquin for type }$C_{n},$ Journal of Algebra, 
\textbf{247}, 295-331 (2002).

\bibitem{Lec3}  \textsc{C.\ Lecouvey,}\textit{\ Combinatorics of crystal
graphs and Kostka-Foulkes polynomials for the root systems }$B_{n},C_{n}$%
\textit{\ and }$D_{n},$ to appear in Journal of European Combinatorics.

\bibitem{Li}  \textsc{D-E. Littlewood,} \textit{The theory of group
characters and matrix representations of groups, }Oxford University Press,
second edition (1958).

\bibitem{Lu}  \textsc{G. Lusztig,} \textit{Singularities, character
formulas, and a }$q$\textit{-analog of weight multiplicities}, Analyse et
topologie sur les espaces singuliers (II-III), Asterisque \textbf{101-102},
208-227 (1983).

\bibitem{mac}  \textsc{I-G. Macdonald,} \textit{Symmetric functions and Hall
polynomials}, Second edition, Oxford Mathematical Monograph, Oxford
University Press, New York, (1995).

\bibitem{NY}  \textsc{A.\ Nakayashiki, Y. Yamada, }\textit{Kostka-Foulkes
polynomials and energy function in sovable lattice models, }Selecta
Mathematica New Series, Vol 3 N$%
%TCIMACRO{\UNICODE[m]{0xb0}}%
%BeginExpansion
{{}^\circ}%
%EndExpansion
$4, 547-599, (1997).

\bibitem{NR}  \textsc{K. Nelsen, A. Ram,} \textit{Kostka-Foulkes polynomials
and Macdonald spherical functions}, preprint (2004), ArXiv: RT/0401298.

\bibitem{OSS}  \textsc{M. Okado, A. Schilling, M. Shimozono, }\textit{A
crystal to rigged configuration bijection for nonexecptional affine
algebras, }Preprint (2002), Arxiv QA/0203163.

\bibitem{ok3}  \textsc{M. Okado, A. Schilling, M. Shimozono, }\textit{%
Virtual Crystals and Fermionic Formulas of Type }$D_{n+1}^{(2)},A_{2n}^{(2)}$%
\textit{\ and }$C_{n}^{(1)},$ Representation Theory, \textbf{7}, 101-163
(2003).

\bibitem{SCS}  \textsc{A. Schilling, M. Shimozono, }$X=M$ for symmetric
powers, Preprint (2004), Arxiv QA/0412376.

\bibitem{SW}  \textsc{A. Schilling, S. O.\ Warnarr, }\textit{Inhomogeneous
lattice paths, generalized Kostka-Foulkes polynomials and }$A_{n-1}$ \textit{%
supernomials, }Comm.\ Math.\ Phys.\ \textbf{202, }359-401 (1999).

\bibitem{sh}  \textsc{M. Shimozono, }\textit{On the }$X=M=K$\textit{\ \
conjecture,} personal communication.

\bibitem{SZ}  \textsc{M. Shimozono, M. Zabrocki, }\textit{Deformed universal
characters for classical and affine algebras, }Preprint\textit{\ }(2004),
ArXiv: CO/0404288.
\end{thebibliography}
\end{document}